\documentclass[red,11pt,a4paper]{article}

\usepackage{amsmath}
\usepackage{amscd}
\usepackage{amssymb}
\usepackage{yfonts}
\usepackage{latexsym}
\usepackage{xcolor}
\usepackage{eufrak}
\numberwithin{equation}{section}
\usepackage[normalem]{ulem}
\usepackage{cases}
\usepackage{comment}
\usepackage{enumitem}
\usepackage{geometry}
\usepackage{cite}
\usepackage[pdfencoding=auto, psdextra]{hyperref}
\usepackage[english]{babel}
\newcommand{\I}{{\mathbb I}}

\newcommand{\CX}{{\mathcal X}}
\newcommand{\CY}{{\mathcal Y}}
\newcommand{\CV}{{\mathcal V}}

\newcommand{\ds}{\,{\rm d}s}

\renewcommand{\div}{{\rm div}\,}

\newcommand{\de}{\partial}

\newcommand{\va}{\textfrak{a}}
\newcommand{\vb}{\textfrak{b}}

\newcommand{\T}{\mathbb{T}}

\newcommand{\N}{\mathbb N}

\newcommand{\Div}{\operatorname{div}}

\newcommand{\linfty}{L_\infty( (0,T) \times \T^3)}
\newcommand{\Linfty}{L_\infty(  \T^3)}

\newcommand{\vw}{{\bf w}}

\newcommand{\vr}{\varrho}

\newcommand{\vt}{\vartheta}

\newcommand{\vu}{\vc{u}}
\newcommand{\vv}{\vc{v}}

\newcommand{\vc}[1]{{\bf #1}}

\newcommand{\Grad}{\nabla}

\newcommand{\pt}{\partial_{t}}

\newcommand{\dx}{{\rm d} {x}}
\newcommand{\dy}{{\rm d} {y}}
\newcommand{\dt}{{\rm d} t }

\newcommand{\lr}[1]{\left( #1 \right)}

\newcommand{\ep}{\varepsilon}

\newcommand{\R}{\mathbb{R}}
\newtheorem{thm}{Theorem}[section]
\newtheorem{lem}[thm]{Lemma}

\title{Regular solutions to the dissipative  Aw-Rascle system}

\author{Nilasis Chaudhuri$^*$ ,\;  Tomasz Piasecki\thanks{Institute of Applied Mathematics and Mechanics, University of Warsaw, ul. Banacha 2, 02-097 Warszawa, Poland}, \; Ewelina Zatorska\thanks{ Mathematics Institute, University of Warwick, Zeeman Building, Coventry CV4 7AL, United Kingdom}
}

\begin{document}
\maketitle 
\begin{abstract}

In this paper we prove the local-in-time existence of regular solutions to dissipative Aw-Rascle system with the offset equal to gradient of some increasing and regular function of density. It is a mixed degenerate parabolic-hyperbolic hydrodynamic model, and we extend the techniques previously developed for compressible Navier-Stokes equations to show the well-posedness of the system in the $L_2-L_2$ setting. We also discuss relevant existence results for offset involving singular or nonlocal  functions of density.
    
\end{abstract}



\section{Introduction}\label{intro}
We investigate the dissipative Aw-Rascle system
\begin{subnumcases}{\label{sys}}
\vr_t + \div(\vr \vu)=0  \\
(\vr \vw)_t + \div(\vr \vw \otimes \vu)=0 
\end{subnumcases}
on $\T^3 \times (0,T)$, where $\T^3 $ is a three-dimensional torus.
The unknown of the system are the density $\vr(t,x)$, and the desired velocity of motion $\vw(t,x)$. The actual velocity of motion $\vu$ is given by the relation:
\begin{equation} \label{uw}
\vu=\vw-\nabla p(\vr),
\end{equation}
where $\Grad p(\vr)$ is the velocity offset, with a given  offset function  $p(\cdot) \in C^5(\R_+)$.

System \eqref{sys} is supplemented with the initial data 
\begin{equation}
    \vr(0,x)=\vr_0(x),\qquad \vu(0,x)=\vu_0(x).
\end{equation}

The purpose of the paper is to prove local-in-time existence of regular solutions to system \eqref{sys} under the following assumptions on the data
\begin{equation} \label{init1}
(\vr_0,\vu_0) \in H^4(\T^3) \times H^3(\T^3).
\end{equation}
We may also assume less regularity of $\vr_0$ at the price of additional assumption on well-prepared data, more precisely
\begin{equation} \label{init2}
\vr_0 \in H^3(\T^3), \; \vu_0+\nabla p(\vr_0) \in H^3(\T^3).    
\end{equation}
%
Our goal is to prove local in time existence of regular solutions to \eqref{sys} 
Note that, using \eqref{uw}, we can rewrite \eqref{sys} as 
\begin{subnumcases}{\label{P:1}}
       \vr_t+\Div(\vr \vw)-\Div(\vr\Grad p(\vr) )=0,\label{P1:1}\\
	(\vr \vw)_t+\Div(\vr \vw\otimes\vw)=\Div(\vr \vw \otimes \nabla p(\vr) ),
	\label{P2:1}
\end{subnumcases} 
which is the equivalent formulation as long as the solution remains sufficiently regular.
Moreover, assuming $\vr>0$ we can further transform this system to obtain 
\begin{subnumcases}{\label{sys3}}
   \vr_t + \div(\vr\vw)=\div(\vr p'(\vr)\nabla\vr),\\
\vw_t + \vu\cdot \nabla \vw = 0, 
\end{subnumcases} 

subject to the initial data
\begin{equation}
    (\vr,\vw)|_{t=0}=(\vr_0,\vu_0+\nabla p(\vr_0)).
\end{equation}

System \eqref{sys} with closure relation \eqref{uw} was recently considered by Acaves et al. \cite{ABDM} in the context of pedestrian flow. Their offset function was actually singular with respect to (w.r.t.) density
\begin{align}\label{singoff}
    p(\vr)= \ep \lr{\frac{1}{\vr}- \frac{1}{\vr_{\rm max}} }^{-\beta}.
\end{align}
This offset function acts as a barrier to ensure that the density remains below its maximum, $\vr_{\rm max}$, which models the formation of congestion within the crowd. 
The same form of offset was previously proposed in the work of Berthelin et al. \cite{BDDR}, as a remedy to the lack of uniform density bounds in the one-dimensional Aw-Rascle road traffic model. For derivation of this model as well as qualitative analysis of the model we refer to\cite{AR, AwKlar}. The classical Aw-Rascle model, however, differs from the system \eqref{sys}, \eqref{uw} not only in its spatial dimension but also because it uses a scalar offset, i.e $ u=w-p(\vr)$. Incorporating the offset as a gradient  in the form of gradient resolves the dimension discrepancy in the closure relation between the velocities $\vu$ and $\vw$. However, the whole system changes its character from hyperbolic to mixed hyperbolic-parabolic type due to additional dissipation effect in the continuity equation of system \eqref{sys3}. The mathematical properties of this system have been explored by Chaudhuri, Gwiazda, and Zatorska in \cite{CGZ}. The authors demonstrated the existence and weak-strong uniqueness of Young-measure solutions to the system  \eqref{sys} \eqref{uw} with $p(\vr)=\vr^\gamma$, $\gamma>0$. Their result states that the measure-valued solution coincides with the strong solution emanating from the same initial data, as long as the latter 
exists.  However, the existence of regular solutions was assumed rather than proven, which motivates the current paper. Our aim is to address this gap.
Initially, we will focus on a generalization of the offset function $p(\vr)=\vr^\gamma$ considered in \cite{CGZ}, followed by an analysis of the well-posedness for two other forms: the singular offset function \eqref{singoff} and a non-local offset function defined as $p(\vr) = K(x) * \vr$. These variations are inspired not only by the aforementioned pedestrian flow model \cite{ABDM}, but also by models that address lubrication effects and collective behaviors, as discussed in \cite{CNPZ, CMPZ, M, CPSZ2024} and related literature.

Lastly, it is important to mention that the well-posedness of the system \eqref{sys}, \eqref{uw} has been previously examined in the framework of weak solutions. Using the method of convex integration it was shown  in \cite{CFZ} that any initial data $(\vr_0,\vu_0)\in C^2(\T^3)\times C^3(\T^3)$ can connect to any terminal data $(\vr_T,\vu_T)\in C^2(\T^3)\times C^3(\T^3)$ consistent with mass and momentum conservation, via a weak solution belonging  the class  
\begin{align*}
    \vr \in C^2([0,T] \times \T^3), \quad \vu \in L^\infty ((0,T) \times \T^3).
\end{align*} 
The corresponding ill-posedness result clearly shows that the existence of so-called wild solutions extends beyond the hyperbolic systems of conservation laws, and  in particular to those experiencing dissipation that degenerates in vacuum.

 In this paper, we extend techniques developed for the compressible Navier-Stokes equations \cite{KNP} to systems of mixed hyperbolic-parabolic type, which exhibit dissipation in the continuity equation but lack it in the momentum equation. We prove the local existence of regular solutions to system (1.1) by applying the method of successive approximations.
 
The main difficulty here is to derive $L_p$ estimates for a linear transport equation. The approach is based on an explicit solution formula obtained by the method of characteristics. Partial results of this type have been used in the theory of compressible Navier-Stokes equations (see among others \cite{VZ},\cite{Zaj},\cite{KNP}), but a consistent $L_p$ theory for transport equations is still missing. Here we address this issue proving quite a general result (Lemma \ref{l:trans}), which may be of independent interest. The dissipativity in \eqref{P2:1} gives parabolic estimates, but a delicate part is to ensure positivity of the solution at each step of the iteration. This issue is addressed in Lemma \ref{lem:CE:exi}.

\subsection{Notation} \label{Notation}
Before stating our main result we shall introduce the notation used in the paper.
\begin{itemize}[leftmargin=*]
    \item Throughout the paper, by $E(\cdot)$ we denote a positive, continuous function such that $E(0)=0$ and $\phi(\cdot)$ denotes a  continuous, positive function. 
    \item We use standard notation $L_p(\T^3)$ and $W^1_p(\T^3)$ for Lebesque and Sobolev spaces on the torus, respectively, and $H^k(\T^3):=W^k_2(\T^3)$. Next, $L_p(0,T;X)$, where $X$ is a Banach space, denotes a Bochner space.
    \item For $T>0$ and $k \in \N$ let us also denote
\begin{equation} \label{def:CX}
\begin{aligned}
&\CX_k(T):=L_2(0,T;H^k(\T^3)) \cap L_\infty(0,T;H^{k-1}(\T^3)),\\
& \CY_k(T):=\{ f\in L_\infty(0,T;H^k(\T^3)): \; f_t \in L_\infty(0,T;H^{k-1}(\T^3))\\
& \CV_k(T):=\{ g \in L_2(0,T;H^{k+1}(\T^3)) \cap L_\infty(0,T;H^k(\T^3)): \;g_t \in L_2(0,T;H^{k-1}(\T^3)) \}
\end{aligned}
\end{equation}
with norms defined in a natural way as appropriate sums of norms.
\end{itemize}

 Since all spaces are considered on the torus, we shall sometimes skip indication of the domain in the definition of space and write $L_p$ instead of $L_p(\T^3)$ etc. .
We are now in a position to state our main result. 
\begin{thm} \label{thm:main}
Assume the initial data satisfies \eqref{init1} or \eqref{init2}. Then there exists $T>0$ such that system \eqref{sys3} admits a unique solution $(\vr,\vw) \in \CV_3(T)\times \CY_3(T)$ with the estimate
$$
\|\vr\|_{\CV_3(T)}+\|\vw\|_{\CY_3(T)} \leq C(\|\vr_0\|_{H^4(\T^3)},\|\vu_0\|_{H^3(\T^3)})
$$
in case of \eqref{init1} or 
$$
\|\vr\|_{\CV_3(T)}+\|\vw\|_{\CY_3(T)} \leq C(\|\vr_0\|_{H^3(\T^3)},\|\vu_0+\nabla p(\vr_0)\|_{H^3(\T^3)})
$$
in case of \eqref{init2}.
\end{thm}

The strategy of the proof involves two main steps:
\begin{itemize}
    \item construction of  solutions to a suitable approximation of system \eqref{sys3},
    \item proof of convergence. 
\end{itemize}
We aim to approximate solutions \eqref{sys3} by solutions to the iterative scheme defined as
\begin{equation} \label{iter}
\left\{ \begin{array}{lr}
\vr^{n+1}_t + \div(\vr^{n+1}\vw^{n+1})=\div(\vr^n p'( \vr^n)\nabla\vr^{n+1}),\\
\vw^{n+1}_t + \vu^n\cdot \nabla \vw^{n+1} = 0,\\
(\vr^{n+1},\vw^{n+1})|_{t=0}= (\vr_0,\vu_0+\nabla p(\vr_0)).
\end{array}\right.
\end{equation}
At each step of the iteration, having $(\vr^n,\vw^n)$ we set $\vu^n=\vw^n-\nabla p(\vr^n)$ and solve the second equation of
\eqref{iter} for $\vw^{n+1}$. Next we use $\vw^{n+1}$ to determine $\vr^{n+1}$ from the first equation. 
Therefore each step of iteration is decoupled to solving separate linear problems   
\begin{equation} \label{lin:diff}
\vr_t + \div(\vr\bar \vw)=\div(\bar \vr p'(\bar \vr)\nabla\vr)\\
\end{equation}
with given $(\bar \vw,\bar \vr)$ and
\begin{equation}\label{lin:trans}
\vw_t + \bar \vu\cdot \nabla \vw = 0
\end{equation}
with given $\bar \vu$. 
Convergence of this iterative scheme is then proved using the Banach fixed point theorem.

\medskip 
The paper is organized as follows.  In Section \eqref{linear} we first solve the linear problems corresponding to the iterative scheme described above in \eqref{lin:diff} and \eqref{lin:trans}. Next, in Section 3, we prove the convergence of the iterative scheme using the contraction argument. Finally, in Section \ref{general-offset} we discuss the existence results for general singular and non-local offset functions; we formulate and prove our other main results -- Theorems \ref{thm:2} and \ref{thm:3}.


\section{Linear theory}\label{linear}
In this section we solve linear problems corresponding to \eqref{lin:diff} and \eqref{lin:trans}.
\subsection{Linear transport equation}  
Consider the linear transport equation 
\begin{equation} \label{trans}
\eta_t + \vv \cdot \nabla \eta =g \;\; {\rm on} \;\; \T^3 \times (0,T) , \quad \eta|_{t=0}=\eta_0 \;\; {\rm on} \;\; \T^3 
\end{equation}
with given vector field $\vv$ and unknown scalar valued $\eta$.
Our goal is to prove the existence of a solution to \eqref{trans} in the regularity framework corresponding to Theorem \ref{thm:main}. We will use the fact that the solution in Lagrangian coordinates 
\begin{equation} \label{lag}
\frac{\de X(t,y)}{\de t}=\vv(t,X(t,y)), \quad X(0,y)=y    
\end{equation}
is constant. The first step is therefore to investigate the regularity properties of solutions to \eqref{lag}. The following 
result improves \cite[Lemma 3.2]{KNP}: 
\begin{lem} \label{l:X}
Assume $\vv \in L_2(0,T;H^3(\T^3))$. Then there exists  a  continuous, positive function $\phi(\cdot)$ denotes such that the solution to \eqref{lag}
satisfies
\begin{align}
&\|\nabla_y X - \I\|_{L_\infty((0,T)\times \T^3)} \leq E(T) \label{x1}\\
&\|\nabla_y X\|_{L_p(0,T;L_\infty(\T^3)}\leq E(T) \quad {\rm for}\;\; 1\leq p<\infty \label{x1b}\\
&\|\nabla_y^2 X\|_{L_\infty(0,T;L_6(\T^3))} \leq \phi(\sqrt{T}\|\vv\|_{L_2(0,T;H^3(\T^3))}) \label{x2}\\
& \|\nabla_y^2 X\|_{L_p(0,T;L_6(\T^3))} \leq E(T) \quad \forall \; 1\leq p<\infty \label{x2b}\\
&\|\nabla_y^3 X\|_{L_\infty(0,T;L_2(\T^3))} \label{x3} \leq \phi(\sqrt{T}\|\vv\|_{L_2(0,T;H^3(\T^3))})\\
&\|\nabla_y^3 X\|_{L_p(0,T;L_2(\T^3))} \leq E(T) \quad \forall \; 1\leq p<\infty \label{x3b}
\end{align}
\end{lem}
{\em Proof.} The first assertion was proved in \cite[Lemma 3.2]{KNP}. Then we have 
$$
\|\nabla_y X\|_{L_p(0,T;L_\infty(\T^3))}\leq T^{1/p}\|\nabla_y X\|_{\linfty}\leq E(T).
$$
In order to prove the bounds for higher derivatives observe that differentiating the solution formula
$$
X(t,y)=y+\int_0^t \vv(s,X(s,y))\ds
$$
with respect to $y$ we obtain 
$$
\nabla_y X(t,y)=y+\int_0^T \nabla_x \vv(x,X(s,y)) \otimes \nabla_y X(s,y)\ds, 
$$
which is equivalent to
$$
\de_t \nabla_y X(t,y) = \nabla_x \vv(x,X(t,y)) \otimes \nabla_y X(t,y).
$$
Differentiating this identity in $y$ we obtain 
\begin{equation} \label{d2x}
\de_t \nabla_y^2 X(t,y) \sim \nabla^2_x \vv(t,X(t,y))(\nabla_y X(t,y))^2 + \nabla_x \vv(t,X(t,y)) \nabla_y^2 X(t,y)
\end{equation}
Multiplying the component corresponding to $\de^2_{y_i y_j} X$
by $|\de^2_{y_i y_j} X|^4 \de^2_{y_i y_j} X$, summing over $i,j$ and integrating over $\T^3$ we get 
\begin{equation} \label{2.1}
\begin{aligned}
\de_t \|\nabla^2_y X(t,\cdot)\|_{L_6(\T^3)}^6 & \leq 
\int_{\T^3} |\nabla^2_x \vv(t,X(t,y))| (\nabla_y X(t,y))^2 |\nabla_y^2 X(t,y)|^5\dy\\
& \int_{\T^3} |\nabla_x \vv(t,X(t,y))| |\nabla_y^2 X(t,y)|^6\dy \\
& \leq \|\nabla^2_x \vv(t,X(t,\cdot))\|_{L_6(\T^3)} \|\nabla_y X(t,y)\|_{L_\infty(\T^3)}^2
\|\nabla_y^2 X(t,\cdot)\|_{L_6(\T^3)}^5 \\
& + \|\nabla_x \vv(t,X(t,\cdot))\|_{L_\infty(\T^3)}\|\nabla_y^2 X(t,\cdot)\|_{L_6(\T^3)}^6
\end{aligned}
\end{equation}
By \eqref{x1}, for small $T$ and any function $f$ of the time variable with values in $L_p(\T^3)$ for $1 \leq p <\infty$ we have 
\begin{equation} \label{eq:norm1}
\begin{aligned}
&\|f(t,X(t,\cdot))\|_{L_p(\T^3)}=
\left( \int_{\T^3}|f(t,X(t,y))|^p \;\dy \right)^{1/p}\\
&=\left( \int_{\T^3}|f(t,X(t,y))|^p |J_y X(t,y)||J_y X(t,y)|^{-1}\;\dy\right)^{1/p}\\ 
&\leq \left(\sup_{y \in \T^3}|J_y X(t,y)|^{-1}\right)^{1/p}
\left(\int_{\T^3}|f(t,x)|^p\;\dx\right)^{1/p}
\leq C \|f(t,\cdot)\|_{L_p(\T^3)},
\end{aligned}
\end{equation}
and similarly 
\begin{equation} \label{eq:norm2}
\|f(t,X(t,\cdot))\|_{L_\infty(\T^3)} \leq C \|f(t,\cdot)\|_{L_\infty(\T^3)}.
\end{equation}
By \eqref{eq:norm1},\eqref{eq:norm2} and Sobolev imbedding, applying \eqref{x1} to the first term of the RHS of \eqref{2.1} we get 
\begin{equation} \label{d2x_1}
\begin{aligned}
\de_t\|\nabla^2_y X(t,\cdot)\|_{L_6(\T^3)}^6 & \leq 
\|\nabla^2_x \vv(t,\cdot)\|_{L_6(\T^3)} \|\nabla_y X(t,y)\|_{L_\infty(\T^3)}^2
\|\nabla_y^2 X(t,\cdot)\|_{L_6(\T^3)}^5 \\
& + \|\nabla_x \vv(t,\cdot)\|_{L_\infty(\T^3)}\|\nabla_y^2 X(t,\cdot)\|_{L_6(\T^3)}^6\\
&\leq C
(\|\nabla^2_x \vv(t,X(t,\cdot))\|_{L_6(\T^3)} + \|\nabla_x \vv(t,X(t,\cdot))\|_{L_\infty(\T^3)})\|\nabla^2_y X(t,\cdot)\|_{L_6(\T^3)}^6\nonumber\\
& \leq C \|\vv(t,X(t,\cdot))\|_{H^3(\T^3)} \|\nabla^2_y X(t,\cdot)\|_{L_6(\T^3)}^6.   
\end{aligned}
\end{equation}
The assumed integrability of $\vv$ allows to conclude \eqref{x2} by Gronwall inequality:
$$
\|\nabla_y^2 X(t,\cdot)\|_{L_6(\T^3)}^6 \leq C \exp\left(\int_0^t\|\vv(s,\cdot)\|_{H^3(\T^3)}\ds \right)
\leq \exp\left(\sqrt{t}\|\vv\|_{L_2(0,t;H^3(\T^3))} \right).
$$
This proves \eqref{x2}, which immedietaly implies \eqref{x2b}.
A remark is due here. In \eqref{d2x_1} we assumed for simplicity that
\begin{equation} \label{Poin}
\|\nabla_y X\|_{L_\infty(\T^3)} \leq C \|\nabla_y^2 X\|_{L_6(\T^3)}
\end{equation}
which does not hold since we don't have Poincar\'e inequality. To make the proof fully precise we would have to replace $\|\nabla_y^2 X\|_{L_6(\T^3)}$ by $\|\nabla_y X\|_{W^1_6(\T^3)}$ which is easy - it is enough to write estimate for $\frac{\de}{\de t}\|\nabla_y X\|_{L_6(\T^3)}^6$. Therefore to avoid additional obvious terms we assume \eqref{Poin}. Similar simplification is also used later in the proof.  

In order to prove \eqref{x3} we differentiate \eqref{d2x} in $y$ obtaining
\begin{equation*} 
\de_t \nabla_y^3 X(t,y) \sim \nabla_x^3\vv(t,X)(\nabla_y X)^3 + \nabla_x\vv(t,X)\nabla_y X\nabla_y^2 X 
+ \nabla_x\vv(t,X)\nabla_y^2 X.
\end{equation*}
Multiplying the equation corresponding to $\de^3_{y_i y_j y_k} X$ by $\de^3_{y_i y_j y_k} X$ and summing over all $i,j,k$ we get 
\begin{equation*}
\begin{aligned}
&\de_t \|\nabla^3_y X(t,\cdot)\|_{L_2(\T^3)}^2 \leq 
\int_{\T^3} \nabla_x^3\vv(t,X(t,y))|\nabla_y X(t,y)|^3|\nabla_y^3 X(t,y)|\,\dy\\
&+\int_{\T^3} |\nabla_x^2 \vv(t,X(t,y))|\; |\nabla_y X(t,y)| |\nabla_y^2 X(t,y)||\nabla_y^3 X(t,y)|\,\dy
+\int_{\T^3}\nabla_x\vv(t,X(t,y))|\nabla_y^3 X(t,y)|^2\,\dy,
\end{aligned}
\end{equation*}
from which, by Sobolev imbedding,\eqref{x1}, \eqref{eq:norm1} and \eqref{eq:norm2}  we obtain 
\begin{equation*}
\begin{aligned}
\de_t \|\nabla^3_y X(t,\cdot)\|_{L_2(\T^3)}^2 
&\leq C \big(\|\nabla_x^3\vv(t,\cdot)\|_{L_2(\T^3)}+\|\nabla_x^2\vv(t,\cdot)\|_{L_6(\T^3)}
+\|\nabla_x \vv(t,\cdot)\|_{L_\infty(\T^3)}\big)\|\nabla^3_y X(t,\cdot)\|_{L_2(\T^3)}^2\\ 
&\leq C \|\vv(t,\cdot)\|_{H^3(\T^3)}\|\nabla^3_y X(t,\cdot)\|_{L_2(\T^3)}^2, 
\end{aligned}
\end{equation*}
and by Gronwall inequality we conclude \eqref{x3}, which implies \eqref{x3b}.
\begin{flushright}
$\square$
\end{flushright}
Now we are in a position to prove a series of estimates for the transport equation \eqref{trans}. As they may be of independent interest, we prove them in possibly general form. 
\begin{lem} \label{l:trans}
Assume $\vv \in \CX_3(T)$ defined in \eqref{def:CX} and $\eta_0 \in H^3(\T^3)$. 
Then the solution to \eqref{trans} satisfies 
\begin{align}
&\|\eta\|_{L_\infty(0,T;W^1_r(\T^3))} \leq \phi(\sqrt{T}\|\vv\|_{L_2(0,T;H^3(\T^3))})\|\eta_0\|_{W^1_r(\T^3)}+E(T)\|g\|_{L_q(0,T;W^1_r(\T^3))}\nonumber\\ 
& \forall\, 1\leq q\leq \infty, \, 1\leq r\leq\infty \label{est:trans:01}\\[5pt]
&\|\nabla_x \eta\|_{L_p(0,T;L_r(\T^3))} \leq E(T)\left( \|\nabla\eta_0\|_{L_\infty(\T^3)}+\|\nabla_x g\|_{L_1(0,T;L_r(\T^3))} \right)\nonumber\\ 
& \forall \, 1\leq p<\infty, \, 1\leq r\leq \infty \label{est:trans:02}\\[5pt]
&\|\nabla_x^2 \eta\|_{L_\infty(0,T;L_6(\T^3))} \leq  \phi(\sqrt{T}\|\vv\|_{L_2(0,T;H^3(\T^3))})\|\eta_0\|_{W^2_6(\T^3)} +E(T)\left( \|\nabla_x g\|_{L_q(0,T;W^1_6(\T^3))} \right) \label{est:trans:03} \\ \nonumber
& \forall\, 1<q\leq \infty,  \\[5pt]
&\|\nabla_x^2 \eta\|_{L_p(0,T;L_6(\T^3))} \leq E(T) \left( \|\eta_0\|_{W^2_6(\T^3)}+ \|\nabla_x g\|_{L_q(0,T;W^1_6(\T^3))} \right) \quad \forall \; 1 \leq p \leq \infty,\label{est:trans:04} \\[5pt] 
&\|\nabla_x^2 \eta\|_{L_\infty(0,T;L_3(\T^3))} \leq  \phi(\sqrt{T}\|\vv\|_{L_2(0,T;H^3(\T^3))})\|\eta_0\|_{W^2_6(\T^3)} +E(T)\left( \|\nabla_x g\|_{L_q(0,T;H^2(\T^3))} \right) \label{est:trans:05} \\ \nonumber
& \forall\, 1<q\leq \infty,  \\[5pt]
&\|\nabla_x^2 \eta\|_{L_p(0,T;L_3(\T^3))} \leq E(T) \left( \|\eta_0\|_{W^2_6(\T^3)}+ \|\nabla_x g\|_{L_q(0,T;W^1_2(\T^3))} \right) \quad \forall \; 1 \leq p < \infty, \label{est:trans:06}\\[5pt] 
&\|\eta\|_{L_\infty(0,T;H^3(\T^3))} \leq \phi(\sqrt{T}\|\vv\|_{\CX_3(T)}) \|\eta_0\|_{W^2_6(\T^3)}\label{est:trans:1}
+E(T) \left( \|\eta_0\|_{H^3(\T^3)}+\|g\|_{L_q(0,T;H^3(\T^3))} \right)
\\ \nonumber 
&\forall\,1<q\leq\infty, 
\\[5pt]
&\|\eta\|_{L_p(0,T;H^3(\T^3))} \leq E(T) \left( \|\eta_0\|_{H^3(\T^3)}+\|g\|_{L_q(0,T;H^3(\T^3))} \right), \label{est:trans:2}\\ \nonumber
&\forall \, 1\leq p<\infty,\, 1<q\leq \infty, \\[5pt] 
&\|\eta_t\|_{L_\infty(0,T;H^2(\T^3))} \leq \phi(\sqrt{T}\|\vv\|_{\CX_3(T)})(\|\eta_0\|_{W^2_6(\T^3)}+\|\nabla g\|_{L_\infty(0,T;H^1(\T^3))}) \nonumber \\
&\hspace{3cm}+E(T)\left( \|\eta_0\|_{H^3(\T^3)}+\|g\|_{L_p(0,T;H^3(\T^3))} \right)  , \label{est:trans:3}\\[5pt]
&\|\eta_t\|_{L_p(0,T;H^2(\T^3))} \leq 
E(T)\left( \|\nabla\eta_0\|_{H^3(\T^3)}+\|g\|_{L_p(0,T;H^3(\T^3))}+\|\nabla g\|_{L_\infty(0,T;H^1(\T^3))} \right). \label{est:trans:4}
\end{align}
\end{lem}

\noindent{\em Proof.}  We have 
$$
\eta(t,X(t,y))=\eta_0(y)+\int_0^t g(s,X(s,y))\ds
$$
Differentiating this identity in $y$ we obtain 
\begin{equation} \label{d1eta}
\nabla_y X(t,y)\nabla_x \eta(t,X(t,y))=\nabla_y \eta_0+\int_0^t \nabla_x g(s,X(s,y))\nabla_y X(s,y) \ds,    
\end{equation}
which, by Lemma \ref{l:X}, implies 
\begin{equation} \label{eta1}
\begin{aligned}
\|\nabla_x \eta\|_{L_\infty(0,T;L_r(\T^3))} & \leq  \|(\nabla_y X)^{-1}\|_{\linfty}\|\nabla_y\eta_0\|_{L_r(\T^3)}\\ 
&+\|(\nabla_y X)^{-1}\int_0^t \nabla_x g(s,X(s,y))\nabla_y X(s,y) \ds \|_{L_\infty(0,T;L_r(\T^3))}\\
&\leq \phi(\sqrt{T}\|\vv\|_{L_2(0,T;H^3(\T^3))})\|\nabla \eta_0\|_{L_r(\T^3)}+E(T)\|\nabla_x g\|_{L_q(0,T;L_r(\T^3))}\\
&\forall \; 1<q\leq\infty, \quad 1\leq r\leq \infty
\end{aligned}
\end{equation}
from which we obtain \eqref{est:trans:01}.
Similarly, using \eqref{x1b} we obtain \eqref{est:trans:02}.
In order to estimate $\nabla_x^2 \eta$ we differentiate \eqref{d1eta} in $y$ to obtain 
\begin{equation*} 
\begin{aligned}
&\frac{\de^2 \eta_0(y)}{\de y_j\de y_k} +\int_0^t \sum_i \left[ g_{x_i}(s,X(s,y))\frac{\de^2 X_i(s,y)}{\de y_j\de y_k}+\left(\sum_l g_{x_ix_l}(s,X(s,y))\frac{\de X_l(s,y)}{\de y_k} \right)\frac{\de X_i(s,y)}{\de y_j}  \right] =\\ 
&=\de_{y_k}\left( \sum_i \frac{\de \eta}{\de x_i}(t,X(t,y))\frac{\de X_i(t,y)}{\de y_j}(t,y) \right)=\\
&=\sum_{i,l}\frac{\de^2 \eta}{\de x_i\de x_l}(t,X(t,y))\frac{\de X_l}{\de y_k}(t,y)\frac{\de X_i}{\de y_j}(t,y)+\nabla_x \eta(t,X(t,y))\cdot \frac{\de^2 X}{\de y_j\de y_k}(t,y)
\end{aligned}
\end{equation*}
for $j,k \in \{1,2,3\}$. 
Rewriting the above system as 
\begin{equation} \label{sys:d2eta}
\begin{aligned}
&\sum_{i,l}\frac{\de^2 \eta}{\de x_i\de x_l}(t,X(t,y))\frac{\de X_l}{\de y_k}(t,y)\frac{\de X_i}{\de y_j}(t,y)=\\
&=\frac{\de^2 \eta_0(y)}{\de y_j\de y_k} 
+\int_0^t \sum_i \left[ g_{x_i}(s,X(s,y))\frac{\de^2 X_i(s,y)}{\de y_j\de y_k}+\left(\sum_l g_{x_ix_l}(s,X(s,y))\frac{\de X_l(s,y)}{\de y_k} \right)\frac{\de X_i}{\de y_j}  \right]\\
&-\nabla_x \eta(t,X(t,y))\cdot \frac{\de^2 X(t,y)}{\de y_j\de y_k}
\end{aligned}
\end{equation}
for $k,j \in \{1,2,3\}$, which is a linear system of $9$ equations for the unknown derivatives 
$\frac{\de^2 \eta}{\de x_i\de x_l}(t,X)$. In order to solve it 
we observe that the diagonal of this system corresponds to $(i,l)=(j,k)$, which means that on the diagonal we have terms $\frac{\de X_k}{\de y_k}\frac{\de X_j}{\de y_j}$, while all entries outside the diagonal contains the terms which are not on the diagonal of $\nabla_y X$. Therefore, by \eqref{x1}, all terms on the diagonal of system \eqref{sys:d2eta}
 are close to $1$ for short times, while all other terms are small. Therefore, system \eqref{sys:d2eta} is uniquely solvable and we obtain  
\begin{equation} \label{eta3aa}
\begin{aligned}
&|\nabla_x^2 \eta(t,X(t,y))| \leq C \Big( |\nabla_y^2 \eta_0| + |\nabla_x \eta(t,X(t,y))||\nabla_y^2 X(t,y)|\\ 
&+ \left| \int_0^t |\nabla_x g(s,X(s,y))||\nabla^2_y X(s,y)|+|\nabla^2_x g(s,X(s,y))||\nabla_y X(s,y)|^2\ds \right| \Big).
\end{aligned}
\end{equation}
By \eqref{x2} and \eqref{eta1} we have 
\begin{equation} \label{eta3a}
\begin{aligned}
&\||\nabla_x\eta| |\nabla^2_y X|\|_{L_\infty(0,T;L_6(\T^3))}\leq 
\|\nabla_x\eta\|_{\linfty}\|\nabla^2_y X\|_{L_\infty(0,T;L_6(\T^3))}\\
&\leq \left( \phi(\sqrt{T}\|\vv\|_{L_2(0,T;H^3(\T^3))})\|\nabla \eta_0\|_{ \Linfty}+E(T)\|\nabla_x g\|_{L_q(0,T;L_\infty(\T^3))} \right)\phi(\sqrt{T}\|\vv\|_{L_2(0,T;H^3(\T^3))})\\
&\leq \phi(\sqrt{T}\|\vv\|_{L_2(0,T;H^3(\T^3))})\|\nabla \eta_0\|_{\Linfty}+E(T)\|\nabla_x g\|_{L_q(0,T;L_\infty(\T^3))} \quad \forall \; 1<q\leq \infty.
\end{aligned}
\end{equation}
Next, by \eqref{x2}
\begin{equation} \label{eta3b}
\begin{aligned}
&\left\|\int_0^t |\nabla_x g| |\nabla^2_y X| \right\|_{L_\infty(0,T;L_6(\T^3))}
\leq \int_0^T \|\nabla_x g\|_{L_\infty(\T^3)} \|\nabla^2_y X\|_{L_6(\T^3)}\\
&\leq \phi(\sqrt{t}\|\vv\|_{L_2(0,T;H^3(\T^3))})\int_0^T \|\nabla_x g\|_{L_\infty(\T^3)}\dt 
\leq E(T) \|\nabla_x g\|_{L_q(0,T;L_\infty(\T^3))} \quad \forall \; 1<q\leq \infty,
\end{aligned}
\end{equation}
and finally 
\begin{equation} \label{eta3c}
\begin{aligned}
&\left\|\int_0^t |\nabla^2_x g| |\nabla_y X|^2\dt \right\|_{L_\infty(0,T;L_p(\T^3))}
\leq \int_0^T \|\nabla^2_x g\|_{L_p(\T^3)}\|\nabla_y X\|_{L_\infty}^2\\ 
&\leq \|\nabla_y X\|_{\linfty}^2 \int_0^T \|\nabla_x^2 g\|_{L_p(\T^3)}\dt 
\leq E(T) \|\nabla^2_x g\|_{L_q(0,T;L_p(\T^3))}\\ 
& \forall \;\; 1<q\leq \infty, \, 1\leq p \leq 6
\end{aligned}
\end{equation}
Combining \eqref{eta3aa},\eqref{eta3a},\eqref{eta3b} and \eqref{eta3c} we obtain \eqref{est:trans:03}.
%
Next, by \eqref{x2b} and \eqref{eta1} we have 
\begin{equation*}
\begin{aligned}
&\left\| |\nabla_x\eta| |\nabla^2_y X| \right\|_{L_p(0,T;L_6(\T^3))}\leq 
\|\nabla_x \eta\|_{\linfty}\|\nabla^2_y X\|_{L_p(0,T;L_6(\T^3))}\\
&\leq E(T) (\textrm{RHS of \eqref{eta1}}) 
\leq E(T) \left( \|\nabla \eta_0\|_{L_\infty(\T^3)}+\|\nabla_x g\|_{L_q(0,T;L_\infty(\T^3))} \right) \quad \forall\; 1<q\leq\infty.
\end{aligned}
\end{equation*}
Combining this estimate with \eqref{eta3aa},\eqref{eta3b} and \eqref{eta3c} 
we arrive at \eqref{est:trans:04}.
Next, similarly to \eqref{eta3a} we obtain
\begin{equation} \label{eta4a}
\begin{aligned}
&\big\|\nabla_x\eta| |\nabla^2_y X|\big\|_{L_\infty(0,T;L_3(\T^3))}\leq 
\|\nabla^2_y X\|_{L_\infty(0,T;L_6(\T^3))} \|\nabla_x\eta\|_{L_1(0,T;L_6(\T^3))}  \\
&\leq \left( \phi(\sqrt{T}\|\vv\|_{L_2(0,T;H^3(\T^3))})\|\nabla \eta_0\|_{\Linfty}+E(T)\|\nabla_x \eta\|_{L_q(0,T;L_6(\T^3))} \right)\phi(\sqrt{T}\|\vv\|_{L_2(0,T;H^3(\T^3))})\\
&\leq \phi(\sqrt{T}\|\vv\|_{L_2(0,T;H^3(\T^3))})\|\nabla \eta_0\|_{\Linfty}+E(T)\|\nabla_x g\|_{L_q(0,T;H^1(\T^3))} \quad \forall \; 1<q\leq \infty.
\end{aligned}
\end{equation}
and, in analogy to \eqref{eta3b},we have
\begin{equation} \label{eta4b}
\begin{aligned}
&\left\|\int_0^t |\nabla_x g| |\nabla^2_y X| \right\|_{L_\infty(0,T;L_3(\T^3))}
\leq \int_0^T \|\nabla_x g\|_{L_6(\T^3)} \|\nabla^2_y X\|_{L_6(\T^3)}\\
&\leq \phi(\sqrt{t}\|\vv\|_{L_2(0,T;H^3(\T^3))})\int_0^T \|\nabla_x g\|_{L_6(\T^3)}\dt 
\leq E(T) \|\nabla_x g\|_{L_q(0,T;H^1(\T^3))} \quad \forall \; 1<q\leq \infty,
\end{aligned}
\end{equation}
Combining \eqref{eta4a}, \eqref{eta4b} and \eqref{eta3c} with $p=2$ we obtain \eqref{est:trans:05}. 
Next, by \eqref{x2b} and \eqref{eta1} we have 
\begin{equation*}
\begin{aligned}
&\left\| |\nabla_x\eta| |\nabla^2_y X| \right\|_{L_p(0,T;L_3(\T^3))}\leq 
\|\nabla_x \eta\|_{L_\infty(L_6(\T^3))}\|\nabla^2_y X\|_{L_p(0,T;L_6(\T^3))}\\
&\leq E(T) (\textrm{RHS of \eqref{eta1}}) 
\leq E(T) \left( \|\nabla \eta_0\|_{\Linfty}+\|\nabla_x g\|_{L_q(0,T;L_6(\T^3))} \right) \quad \forall\; 1<q\leq\infty,
\end{aligned}
\end{equation*}
which combined with \eqref{eta4b} and \eqref{eta3c} for $p=2$ gives \eqref{est:trans:06}.

In order to estimate the third order derivatives we differentiate \eqref{sys:d2eta} w.r.t $y_m$, which yields 
\begin{equation} \label{sys:d3eta}
\begin{aligned}
&\sum_{i,l,n}\frac{\de^2\eta}{\de x_i\de x_l\de x_n}(t,X)\frac{\de X_n}{\de y_m}\frac{\de X_l}{\de y_k}\frac{\de X_i}{\de y_j}=\\
&=\frac{\de^3 \eta_0}{\de y_j\de y_k\de y_m}
-\de_{y_m}\left( \nabla_x \eta(t,X)\frac{\de^2 X}{\de y_j\de y_k}(t,y) \right)\\
&-\sum_{i,l}\frac{\de^2 \eta}{\de x_i \de x_l}(t,X)\left( \frac{\de^2 X_l}{\de y_k \de y_m}\frac{\de X_i}{\de y_l}+\frac{\de X_i}{\de y_l\de y_m}\frac{\de X_l}{\de y_k} \right)(t,y)\\
&+\int_0^t \de_{y_m} \left\{ \sum_i \left[ g_{x_i}(s,X)\frac{\de^2 X_i}{\de y_j\de y_k}+\left(\sum_l g_{x_ix_l}(s,X)\frac{\de X_l}{\de y_k} \right)\frac{\de X_i}{\de y_j}  \right] \right\},
\end{aligned}
\end{equation} 
where $X=X(t,y)$ or $X=X(s,y)$ according to \eqref{sys:d2eta}. 
Similarly as in case of \eqref{sys:d2eta}, it is a system of $27$ linear equations for the third order derivatives of $\eta$. On the diagonal we have terms corresponding to $(i,l,n)=(j,k,m)$, which, again by \eqref{x1}, are close to one, while all other entries are small for small times. Therefore \eqref{sys:d3eta} is uniquely solvable and we obtain 
\begin{equation} \label{eta5:0}
\begin{aligned}
|\nabla^3_x \eta|\leq & C \left( |\nabla^3_y \eta_0|+|\nabla^2_x \eta| |\nabla_y X| |\nabla_y^2 X| + |\nabla_x \eta| |\nabla^3_y X| \right. \\
&\left. +\int_0^t |\nabla_x g||\nabla_y^3 X|+|\nabla^2_x g||\nabla_y X||\nabla^2_y X|+|\nabla_x^3 g||\nabla_y X|^3\dt
\right)
\end{aligned}
\end{equation}
Let us estimate the RHS of \eqref{eta5:0}. For the second term, by \eqref{x2} and \eqref{est:trans:03} we have 
\begin{equation} \label{eta5:2}
\begin{aligned}
&\left\| |\nabla^2_x \eta| |\nabla_y X| |\nabla^2_y X| \right\|_{L_\infty(0,T;L_2(\T^3))}
\leq C \|\nabla^2_x \eta\|_{L_\infty(0,T;L_3(\T^3))}  \|\nabla^2_y X\|_{L_\infty(0,T;L_6(\T^3))}\\[5pt]
&\leq \phi(\sqrt{t}\|\vv\|_{L_2(0,T;H^3(\T^3))})\,\left[\textrm{RHS of \eqref{est:trans:03}}\right]\\[5pt]
&\leq \phi(\sqrt{t}\|\vv\|_{L_2(0,T;H^3(\T^3))})\|\eta_0\|_{W^2_6(\T^3)}
+E(T)\left( \|\nabla_x g\|_{L_q(0,T;L_\infty(\T^3))} + \|\nabla^2_x g\|_{L_q(0,T;L_6(\T^3))} \right),
\end{aligned}
\end{equation}
and for the third, by \eqref{x3} and \eqref{est:trans:01}  
\begin{equation*} \label{eta5:3}
\begin{aligned}
&\left\| |\nabla_x\eta| |\nabla^3_y X|  \right\|_{L_\infty(0,T;L_2(\T^3))}
\leq \|\nabla_x\eta\|_{\linfty}  \|\nabla^3_y X\|_{L_\infty(0,T;L_2(\T^3))}\\[5pt]
& \phi(\sqrt{T}\|\vv\|_{L_2(0,T;H^3(\T^3))})\,\left[\textrm{RHS of \eqref{est:trans:01}}\right]\\[5pt]
&\leq \phi(\sqrt{T}\|\vv\|_{L_2(0,T;H^3(\T^3))})\|\eta_0\|_{W^1_\infty}+E(T)\|g\|_{L_q(0,T;W^1_\infty(\T^3))}. 
\end{aligned}
\end{equation*}
It remains to estimate the terms with $g$. By \eqref{x3} we have 
\begin{equation} \label{eta5:4}
\begin{aligned}
&\left\| \int_0^t |\nabla_x g||\nabla_y^3 X| \right\|_{L_\infty(0,T;L_2(\T^3))} \leq \int_0^T \|\nabla_x g(t,\cdot)\|_{L_\infty}\|\nabla_y^3 X(t,\cdot)\|_{L_2(\T^3)}\dt\\[5pt]
&\leq E(T)\|\nabla_x g\|_{L_q(0,T;L_\infty(\T^3))} \quad \forall\, q<1\leq\infty.
\end{aligned}
\end{equation}
Next, by \eqref{x2},
\begin{equation} \label{eta5:5}
\begin{aligned}
&\left\| \int_0^t|\nabla^2_x g||\nabla_y X||\nabla^2_y X| \right\|_{L_\infty(0,T;L_2(\T^3))} 
\leq C \int_0^T \|\nabla^2_x g(t,\cdot)\|_{L_3(\T^3)}\|\nabla^2_y X(t,\cdot)\|_{L_6(\T^3)}\\[5pt]
&\leq C \|\nabla^2_y X\|_{L_{\infty}(0,T;L_6(\T^3))} \int_0^T \|\nabla^2_x g(t,\cdot)\|_{L_3(\T^3)}\dt
\leq E(T) \|\nabla^2_x g\|_{L_q(0,T;L_3(\T^3))} \quad \forall\, 1<q\leq\infty,
\end{aligned}
\end{equation}
and finally 
\begin{equation} \label{eta5:6}
\begin{aligned}
\left\| |\nabla^3_x g||\nabla_y X|^3| \right\|_{L_\infty(0,T;L_2(\T^3))}
\leq C\int_0^T \|\nabla_x^3 g(t,\cdot)\|_{L_2(\T^3)}\dt \leq E(T)\|\nabla_x^3 g\|_{L_q(0,T;L_2(\T^3))} \quad \forall\, 1<q\leq\infty.
\end{aligned}
\end{equation}
Combining \eqref{eta5:0}-\eqref{eta5:6} and applying Sobolev imbedding to estimate all terms containing $g$ by a single norm we obtain 
\begin{equation*} \label{eta5}
\begin{aligned}
&\|\nabla_x^3 \eta\|_{L_\infty(0,T;L_2(\T^3))} \leq \phi(\sqrt{T}\|\vv\|_{L_2(0,T;H^3(\T^3))}) \|\eta_0\|_{W^2_6(\T^3)}\\[5pt] 
&+E(T) \left( \|g\|_{L_q(0,T;H^3(\T^3))} \right) \quad \forall\,1<q\leq\infty,
\end{aligned}
\end{equation*}
which, together with estimates on lower order derivatives of $\eta$ completes the proof of \eqref{est:trans:1}.

In order to prove of \eqref{est:trans:2} observe that, for any finite $p$, by \eqref{x2} and \eqref{est:trans:04} we have 
\begin{equation} \label{eta6:1}
\begin{aligned}
\|\nabla_x^2 \eta\nabla_y X\nabla_y^2 X\|_{L_p(0,T;L_2(\T^3))}  
\leq C \|\nabla_y^2 X\|_{L_\infty(0,T;L_6(\T^3))} \|\nabla_x^2 \eta(t)\|_{L_p(0,T;L_3(\T^3))} \\[5pt]
\leq E(T) \left( \|\eta_0\|_{W^2_6(\T^3)}+ \|\nabla_x g\|_{L_q(0,T;L_\infty)}+\|\nabla^2_x g\|_{L_q(0,T;L_6(\T^3))} \right) \quad \forall 1<q\leq \infty.
\end{aligned}
\end{equation}
Similarly by  \eqref{x3} and \eqref{est:trans:2} we obtain 
\begin{equation} \label{eta6:2}
\begin{aligned}
&\| \nabla_x \eta \nabla_y^3 X \|_{L_p(0,T;L_2(\T^3))} \leq C \|\nabla_y^3 X\|_{L_\infty(0,T;L_2(\T^3))} \|\nabla_x \eta\|_{L_p(0,T;L_\infty(\T^3))}\\ 
&\leq E(T) \left( \|\nabla\eta_0\|_{\Linfty}+\|\nabla_x g\|_{L_1(0,T;L_\infty(\T^3))}\right).
\end{aligned}
\end{equation}
For the terms with $g$ on the RHS of \eqref{eta5:0} we use the estimates \eqref{eta5:4}-\eqref{eta5:6}. Combining them with \eqref{eta6:1}-\eqref{eta6:2} we obtain 
\begin{equation*}
\|\nabla_x^3 \eta\|_{L_p(0,T;L_2(\T^3))} \leq E(T) \left( \|\eta_0\|_{H^3(\T^3)}+\|g\|_{L_q(0,T;H^3(\T^3))} \right) 
\quad \forall \, 1\leq p<\infty,\, 1<q\leq \infty,     
\end{equation*}
which completes the proof of \eqref{est:trans:2}. 
Now we can use \eqref{trans} to prove the estimates for $\eta_t$. First we immediately get 
\begin{equation*} 
\|\eta_t\|_{\linfty} \leq \phi(\sqrt{t}\|\vv\|_{\CX_3(T)}), \quad 
\|\eta_t\|_{L_p(0,T;L_\infty(\T^3))} \leq E(T) \quad \forall \; 1\leq p<\infty.
\end{equation*}
Next we differentiate \eqref{trans} in the space variable to obtain
\begin{equation} \label{d1trans}
\nabla \eta_t \sim \nabla \vv \nabla \eta + \vv \nabla^2 \eta+\nabla g. 
\end{equation}
By \eqref{est:trans:01} we have 
\begin{equation} \label{eta7:1a}
\begin{aligned}
&\|\nabla \vv \nabla \eta\|_{L_\infty(0,T;L_6(\T^3))}\leq \|\nabla \vv\|_{L_\infty(0,T;L_6(\T^3))}\|\nabla \eta\|_{\linfty}\\
&\leq\phi(\sqrt{T}\|\vv\|_{L_2(0,T;H^3(\T^3))})\|\eta_0\|_{W^1_\infty(\T^3)}+E(T)\|g\|_{L_q(0,T;W^1_\infty(\T^3))} \quad \forall\, 1<q\leq \infty,
\end{aligned}
\end{equation}
by \eqref{est:trans:03}
\begin{equation} \label{eta7:2a}
\begin{aligned}
&\|\vv \nabla^2 \eta\|_{L_\infty(0,T;L_6(\T^3))}\leq \| \vv\|_{\linfty}\|\nabla^2 \eta\|_{L_\infty(0,T;L_6(\T^3))} \\ 
&\leq \phi(\sqrt{T}\|\vv\|_{L_2(0,T;H^3(\T^3))})\|\eta_0\|_{W^2_6(\T^3)}
+E(T)\left( \|\nabla_x g\|_{L_q(0,T;L_\infty(\T^3))}+\|\nabla^2_x g\|_{L_q(0,T;L_6(\T^3))} \right),   \\ 
\quad \forall\, 1<q\leq \infty,
\end{aligned}
\end{equation}
by \eqref{est:trans:02}
\begin{equation} \label{eta7:1b}
\begin{aligned}
&\|\nabla \vv\nabla\eta(t)\|_{L_p(0,T;L_6(\T^3))} \leq \|\nabla \vv\|_{\linfty} \|\nabla\eta(t)\|_{L_p(0,T;L_6(\T^3))}\\
&\leq E(T)\left( \|\nabla\eta_0\|_{L_\infty(\T^3)}+\|\nabla_x g\|_{L_1(0,T;L_\infty(\T^3))} \right),
\end{aligned}
\end{equation}
and finally by \eqref{est:trans:04}
\begin{equation} \label{eta7:2b}
\begin{aligned}
&\int_0^T \|\vv(t,\cdot)\nabla^2 \eta(t,\cdot)\|_{L_6(\T^3)}^p\;\dt \leq \|\vv\|_{\linfty}^p \int_0^T \|\nabla^2 \eta(t,\cdot)\|_{L_6(\T^3)}^p\;\dt\\ 
&\leq E(T) \left( \|\eta_0\|_{W^2_6(\T^3)}+ \|\nabla_x g\|_{L_q(0,T;W^1_6(\T^3))} \right) \quad \forall\, 1\leq p<\infty, \,1<q\leq\infty, 
\end{aligned}
\end{equation}
so alltogether we obtain 
\begin{equation*} 
\begin{aligned}
&\|\nabla \eta_t\|_{L_\infty(0,T;L_6(\T^3))} \leq \phi(\|\vv\|_{\CX_3(T)})(\|\eta_0\|_{W^2_6(\T^3)}+\|\nabla_x g\|_{L_\infty(0,T;L_6(\T^3))})+E(T)\|\nabla_x g\|_{L_q(0,T;W^1_6(\T^3))} \\ 
& \quad \forall 1<q\leq\infty,\\[5pt]
&\|\nabla \eta_t\|_{L_p(0,T;L_6(\T^3))} \leq E(T)\left( \|\eta_0\|_{W^2_6(\T^3)}+\|\nabla g\|_{L_q(0,T;W^1_6(\T^3))}+\|\nabla g\|_{L_\infty(0,T;L_6(\T^3))} \right) \\
& \quad \forall 1\leq p<\infty,\,1<q\leq\infty.
\end{aligned}
\end{equation*}
Finally we differentiate \eqref{d1trans} once more in space:
$$
\nabla^2\eta_t \sim \nabla^2 \vv \nabla\eta + \nabla\vv\nabla^2\eta+\vv\nabla^3 \eta+\nabla^2 g.
$$
For the first term we have, by \eqref{est:trans:01}, 
\begin{equation} \label{eta8:1a}
\begin{aligned}
&\|\nabla^2 \vv \nabla\eta\|_{L_\infty(0,T;L_2(\T^3))} \leq \|\nabla \eta\|_{\linfty} \|\nabla^2 \vv\|_{L_\infty(0,T;L_2(\T^3))}\\
&\leq \phi(\sqrt{T}\|\vv\|_{L_2(0,T;H^3(\T^3))})\|\eta_0\|_{W^1_\infty(\T^3)}+E(T)\|g\|_{L_q(0,T;W^1_\infty(\T^3))} \quad \forall\, 1<q\leq \infty, 
\end{aligned}
\end{equation}
and, by \eqref{est:trans:02},
\begin{equation} \label{eta8:1b}
\begin{aligned}
&\|\nabla^2 \vv \nabla \eta(t)\|_{L_p(0,T;L_2(\T^3))}^p \leq \|\nabla \eta\|_{L_p(0,T;L_\infty(\T^3))} \int_0^T \|\nabla^2 \vv(t)\|_{L_\infty(0,T;L_2(\T^3))}\\
&\leq E(T)\left( \|\nabla\eta_0\|_{L_\infty(\T^3)}+\|\nabla_x g\|_{L_1(0,T;L_\infty(\T^3))} \right) \quad \forall \; 1 \leq p<\infty.
\end{aligned}
\end{equation}
For the second term, by \eqref{est:trans:03}, 
\begin{equation} \label{eta8:2a}
\begin{aligned}
&\|\nabla\vv\nabla^2\eta\|_{L_\infty(0,T;L_2(\T^3))}\leq \|\nabla\vv\|_{L_\infty(0,T;L_6(\T^3))} \|\nabla^2\eta\|_{L_\infty(0,T;L_3(\T^3))}\\
&\leq \phi(\sqrt{T}\|\vv\|_{L_2(0,T;H^3(\T^3))})\|\eta_0\|_{W^2_6(\T^3)} +E(T)\left( \|\nabla_x g\|_{L_q(0,T;W^1_6(\T^3))} \right) \quad \forall\, 1<q\leq\infty
\end{aligned}
\end{equation}
and, by \eqref{est:trans:04}
\begin{equation} \label{eta8:2b}
\begin{aligned}
&\|\nabla \vv \nabla^2 \eta(t)\|_{L_p(0,T;L_2(\T^3))} \leq \|\nabla \vv\|_{L_\infty(0,T;L_6(\T^3))} \|\nabla^2 \eta(t)\|_{L_p(0,T;L_6(\T^3))}\\ &\leq E(T) \left( \|\eta_0\|_{W^2_6(\T^3)}+ \|\nabla_x g\|_{L_q(0,T;W^1_6(\T^3))} \right)  \quad
\forall \; 1 \leq p <\infty,\, 1 < q \leq \infty.
\end{aligned}
\end{equation}
Finally, to estimate the last term we apply \eqref{est:trans:1} to get 
\begin{equation} \label{eta8:3a}
\begin{aligned}
&\|\vv\nabla^3\eta\|_{L_\infty(0,T;L_2(\T^3))}\leq \|\vv\|_{\linfty} \|\nabla^3\eta\|_{L_\infty(0,T;L_2(\T^3))}\\
&\leq \phi(\sqrt{t}\|\vv\|_{\CX_3(T)}) \|\eta_0\|_{W^2_6(\T^3)}
+E(T) \left( \|\eta_0\|_{H^3(\T^3)}+\|g\|_{L_q(0,T;H^3(\T^3))} \right)\quad\forall 1<q\leq\infty
\end{aligned}
\end{equation}
and \eqref{est:trans:2} to obtain
\begin{equation} \label{eta8:3b}
\begin{aligned}
&\|\vv \nabla^3 \vr(t)\|_{L_p(0,T;L_2(\T^3))} \leq \| \vv\|_{\linfty}^p \int_0^T \|\nabla^3 \eta(t)\|_{L_p(0,T;L_2(\T^3))}\\
&\leq E(T) \left( \|\eta_0\|_{H^3(\T^3)}+\|g\|_{L_q(0,T;H^3(\T^3))} \right)\quad
\forall \, 1\leq p<\infty,\, 1<q\leq \infty.
\end{aligned}
\end{equation}
Combining \eqref{eta7:1a},\eqref{eta7:2a},\eqref{eta8:1a},\eqref{eta8:2a} and \eqref{eta8:3a} we obtain \eqref{est:trans:3}. Finally, \eqref{eta7:1b},\eqref{eta7:2b},\eqref{eta8:1b},\eqref{eta8:2b} and \eqref{eta8:3b} allows to conclude \eqref{est:trans:4}, which completes the proof.  
\begin{flushright}
$\square$
\end{flushright}

\subsection{Linear continuity equation with dissipation} 
In this section we investigate the linear problem 
\begin{equation} \label{CE:lin}
\vr_t + \div(\vr \vv)-\div(\va\nabla\vr)=\vb, \quad \vr|_{t=0}=\vr_0
\end{equation}
Concerning the regularity of the data, we keep in mind that 
the above system corresponds to the first equation of \eqref{sys3}. Therefore, taking into account Lemma \ref{l:trans},
it is sufficient to assume $\vv \in L_\infty(0,T;H^3(\T^3))$. The parabolic maximal regularity then leads to the following result 
\begin{lem}\label{lem:CE:exi}
Assume $\vr_0 \in H^2(\T^3)$, $\vv \in L_\infty(0,T;H^2(\T^3))$, {$\va \in L_\infty(0,T;W^1_\infty(\T^3)) \cap L_2(0,T;H^2(\T^3))$} with $\va \geq c>0$, and $\vb \in L_2(0,T;H^1(\T^3))$. Then \eqref{CE:lin} admits a unique solution satisfying 
\begin{equation} \label{est:CE:0}
\begin{aligned}
&\|\vr\|_{\CV_2(T)} \leq C \Big[ \|\vr_0\|_{H^2(\T^3)} +\|\vb\|_{L_2(0,T;H^1(\T^3))}\\ 
&+\|\vr\|_{L_\infty(0,T;H^2(\T^3))} 
\Big( T \big(\|\vv\|_{L_\infty(0,T;H^2(\T^3))}
+ \|\nabla\va\|_{\linfty}\big) +\|\nabla^2 \va\|_{L_2(0,T;L_2(\T^3))}]\Big)\Big].
\end{aligned}
\end{equation}
If $\vr_0 \in H^3(\T^3)$, $\vb \in L_2(0,T;H^2(\T^3))$  $\vv \in L_\infty(0,T;H^3(\T^3))$ and $\va \in L_\infty(W^1_\infty(\T^3)) \cap L_2(0,T;H^3(\T^3))$ with $\va \geq c>0$, then  
\begin{equation} \label{est:CE}
\begin{aligned}
&\|\vr\|_{\CV_3(T)} \leq C \Big[ \|\vr_0\|_{H^3(\T^3)} +\|\vb\|_{L_2(0,T;H^2(\T^3))}\\ 
&+\|\vr\|_{L_\infty(0,T;H^3(\T^3))} 
\Big( T \big(\|\vv\|_{L_\infty(0,T;H^3(\T^3))}
+ \|\nabla\va\|_{\linfty}\big) +\|\nabla^2 \va\|_{L_2(0,T;H^1(\T^3))}]\Big)\Big].
\end{aligned}
\end{equation}
Moreover, if we assume 
\begin{equation} \label{b:form}
\vb=\div\lr{\vr \mathbf{b} } \;\; {\rm with} \;\; \Vert \div \mathbf{b}\Vert_{L_1(0,T;L_\infty(\T^3))} \leq C 
\end{equation} 
and the initial data is strictly positive, i.e., $\inf_{x\in \T^3} \vr_0 >0$ then we have
\begin{equation} \label{min:vr}
\min_{(t,x)\in [0,T]\times \T^3 }\vr(t,x) >0 . 
\end{equation}
\end{lem}
{\em Proof.}  Rewriting \eqref{CE:lin} as 
\begin{equation} \label{para:1}
\vr_t - \div(\va\nabla\vr) = -\div (\vr\vv)+\vb    
\end{equation}
we immediately obtain the bound 
\begin{equation} \label{est:1}
\begin{aligned}
&\|\vr\|_{L_2(0,T;H^2(\T^3))}+\|\vr\|_{L_\infty(0,T;H^1(\T^3))}+\|\vr_t\|_{L_2(0,T;L_2(\T^3))}\\ 
&\leq C\left(\|\vr_0\|_{H^1(\T^3)}+ \|\div(\vr\vv)\|_{L_2(0,T;L_2(\T^3))}+\|\vb\|_{L_2(0,T;L_2(\T^3))}\right).
\end{aligned}
\end{equation}
Differentiating \eqref{para:1} in $x_i$ we obtain 
\begin{equation*}
(\vr_{x_i})_t - \div(\va \nabla\vr_{x_i})=-(\div(\vr\vv))_{x_i} + \div(\va_{x_i}\nabla \vr)+\vb_{x_i} =: F_1^i,    
\end{equation*}
therefore 
\begin{equation} \label{est:2}
\begin{aligned}
&\|\vr_{x_i}\|_{L_2(0,T;H^2(\T^3))}+\|\vr_{x_i}\|_{L_\infty(0,T;H^1(\T^3))}+\|\de_t\vr_{x_i}\|_{L_2(0,T;L_2(\T^3))}\\ 
&\leq C\left(\|\vr_0\|_{H^2(\T^3)}+\|F_1^i\|_{L_2(0,T;L_2(\T^3))}\right).
\end{aligned}
\end{equation}
Under the assumed regularity of $\vv$ and $\va$, using \eqref{est:1} we can find appropriate bound on $\|F_1\|_{L_2(0,T;L_2(\T^3))}$. Namely, 
$$
F_1 \sim \vv\nabla^2\vr + \nabla \vv\nabla\vr+\vr\nabla^2\vv+\nabla^2\va\nabla\vr+\nabla\va\nabla^2\vr+\nabla \vb.
$$
We have
\begin{align*}
&\|\vv\nabla^2\vr\|_{L_2(0,T,L_2(\T^3))}\leq T \|\vv\|_{\linfty}\|\nabla^2\vr\|_{L_\infty(0,T;L_2(\T^3))},\\[3pt]
&\|\nabla\vv\nabla\vr\|_{L_2(0,T;L_2(\T^3))}\leq T\|\nabla\vv\|_{\linfty}\|\nabla\vr\|_{L_\infty(0,T;L_2(\T^3))},\\[3pt]
&\|\vr\nabla^2\vv\|_{L_2(0,T;L_2(\T^3))}\leq T\|\vr\|_{\linfty}\|\nabla^2\vv\|_{L_\infty(0,T;L_2(\T^3))},\\[3pt]
& \|\nabla \vr\nabla^2\va\|_{L_2(0,T;L_2(\T^3))}\leq \|\vr\|_{\linfty}\|\nabla^2\va\|_{L_2(0,T;L_2(\T^3))},\\[3pt]
& \|\nabla^2 \vr\nabla\va\|_{L_2(0,T;L_2(\T^3))}\leq T\|\nabla\va\|_{\linfty}\|\nabla^2\vr\|_{L_\infty(0,T;L_2(\T^3))},
\end{align*}
which gives 
\begin{equation*}
\begin{aligned}
\|F_1\|_{L_2(0,T;L_2(\T^3))}\leq \|\vr\|_{L_\infty(0,T;H^2(\T^3))}\Big(& T[\|\vv\|_{L_\infty(0,T;H^2(\T^3))}+\|\nabla\va\|_{\linfty}]\\
&+\|\nabla^2\va\|_{L_2(0,T;L_2(\T^3))}+\|\vb\|_{L_2(0,T;H^1(\T^3))} \Big),
\end{aligned}
\end{equation*}
which together with \eqref{est:2} implies \eqref{est:CE:0}.

Next, let $\alpha$ be any multiindex with $|\alpha|=2$. Applying $D^\alpha$ to \eqref{para:1} we obtain 
\begin{equation} \label{para:2}
(D^\alpha \vr)_t - \div (\va\nabla D^\alpha\vr) = F_2^\alpha,
\quad D^\alpha \vr|_{t=0}= D^\alpha \vr_0,
\end{equation}
where 
$$
F_2^\alpha \sim \vv \nabla^3 \vr + \nabla \vv\nabla^2\vr+\nabla^2\vv\nabla\vr+\vr\nabla^3\vv
+\sum_{k=1}^3\nabla^k\va \nabla^{4-k}\vr+\nabla^2 \vb.
$$
In order to prove \eqref{est:CE}
we have to estimate the $L_2(0,T;L_2(\T^3))$ norm of $F_2^\alpha$. We have 
\begin{align*}
&\|\vv\nabla^3\vr\|_{L_2(0,T;L_2(\T^3))} 
\leq T \|\vv\|_{\linfty}\|\nabla^3\vr\|_{L_\infty(0,T;L_2(\T^3))},\\
& \|\nabla \vv\nabla^2\vr\|_{L_2(0,T;L_2(\T^3))} \leq T \|\nabla \vv\|_{\linfty}\|\nabla^2\vr\|_{L_\infty(0,T;L_2(\T^3))},\\
& \|\nabla^2 \vv\nabla \vr\|_{L_2(0,T;L_2(\T^3))} 
\leq T \|\nabla^2\vv\|_{L_\infty(0,T;L_6(\T^3))}\|\nabla\vr\|_{L_\infty(0,T;H_2(\T^3))},\\
&\|\vr\nabla^3\vv\|_{L_2(0,T;L_2(\T^3))} \leq T \|\vr\|_{\linfty}\|\nabla^3\vv\|_{L_\infty(0,T;L_2(\T^3))}.
\end{align*}
Combining these estimates with Sobolev imbedding we obtain 
\begin{equation} \label{est:r1}
\|\vv \nabla^3 \vr + \nabla \vv\nabla^2\vr+\nabla^2\vv\nabla\vr+\vr\nabla^3\vv\|_{L_2(0,T;L_2(\T^3))} \leq T \|\vr\|_{L_\infty(0,T;H^3(\T^3))}\|\vv\|_{L_\infty(0,T;H^3(\T^3))}   
\end{equation}
The terms with $\va$ can be treated as follows 
\begin{align*}
&\|\nabla\vr\nabla^3\va\|_{L_2(0,T;L_2(\T^3))}^2 \leq \|\nabla \vr\|^2_{\linfty} \|\nabla^3\va\|_{L_2(0,T;L_2(\T^3))}^2,\\
&\|\nabla^2\vr\nabla^2\va\|_{L_2(0,T;L_2(\T^3))}^2 \leq 
\|\nabla^2\vr\|_{L_\infty(0,T;L_6(\T^3))}^2 \|\nabla^2\va\|_{L_2(0,T;L_3(\T^3))}^2,\\
&\|\nabla^3\vr \nabla \va\|_{L_2(0,T;L_2(\T^3))}^2 \leq 
T \|\nabla\va\|^2_{\linfty}\|\nabla^3 \vr\|_{L_\infty(0,T;L_2(\T^3))},
\end{align*}
which together with Sobolev imbedding yields 
\begin{equation} \label{est:r2}
\begin{aligned}
&\|\sum_{k=1}^3\nabla^m\va \nabla^{4-m}\vr\|_{L_2(0,T;L_2(\T^3))}\\
&\leq \|\vr\|_{L_\infty(0,T;H^3(\T^3))} \big[T\|\nabla\va\|_{\linfty}+\|\nabla^2\va\|_{L_2(0,T;H^1(\T^3))}\big].
\end{aligned}
\end{equation}
Combining \eqref{est:r1} and \eqref{est:r2} with the maximal regularity estimate for \eqref{para:2}
we obtain
\begin{equation*}
\begin{aligned}
&\|\de_t\nabla^2 \vr\|+\|\nabla^2\vr\|_{L_\infty(0,T;H^1(\T^3))}+\|\nabla^2\vr\|_{L_2(0,T;H^2(\T^3))} 
\leq C \Big[\|\vr_0\|_{H^3(\T^3)}+\|\vb\|_{L_2(0,T;H^2(\T^3))} \\[5pt] 
&+\|\vr\|_{L_\infty(0,T;H^3(\T^3))} \Big( T \big(\|\vv\|_{L_\infty(0,T;H^3(\T^3))}
+ \|\nabla\va\|_{\linfty}\big) +\|\nabla^2\va\|_{L_2(0,T;H^1(\T^3))}\Big)\Big].
\end{aligned}
\end{equation*}
Combining this estimate with
\eqref{est:1} and \eqref{est:2} we obtain \eqref{est:CE}.

It remains to prove \eqref{min:vr} under additional assumption \eqref{b:form}.
Consider a function $\psi : \R \rightarrow \R$ as 
\[\mathit{\psi(\lambda)}= 
   \begin{cases}
   \frac{1}{2}\lambda^2, & \lambda\leq 0\\
   0, &\lambda> 0.\\
\end{cases}\]
   Clearly, $\psi \in C^{1}(\R)$ with
\[\mathit{\psi^{'}(\lambda)}= 
   \begin{cases}
   \lambda, &\lambda\leq 0\\
   0, &\lambda> 0.\\
\end{cases}\]
Let us denote $K= \underline{\vr_0}:= \min_{x \in \T^3} \vr_0 >0 $.
Now, we consider a function $h_z:[0,T)\rightarrow \R$ as
\begin{align*}
h_z(t)=\int_{\T^3}\psi(\vr(t,x)-K)\; \dx. 
\end{align*}   
Since,  $\vr \in \CV_3(T) $, a straightforward computation yields
\begin{align*}
    h_z^{\prime}(t)= \int_{\T^3}\psi^{\prime}(\vr(t,x)-K) \pt \vr(t,x) \; \dx.
\end{align*}
Using the equation \eqref{CE:lin}, assumption \eqref{b:form} and integration by parts, we obtain 
\begin{align*}
    h_z^{\prime}(t)= \int_{\T^3}\psi^{\prime \prime}(\vr(t,x)-K) \vr \lr{\vv-\mathbf{b}} \nabla(\vr(t,x)-K) \; \dx - \int_{\T^3} \va \psi^{\prime \prime}(\vr(t,x)-K) \vert \nabla_x \vr(t,x) \vert^2  \; \dx .
\end{align*}
From the assumption on $\va$, we get 
\[\int_{\T^3} \va \psi^{\prime \prime}(\vr(t,x)-K) \vert \nabla (\vr(t,x)) \vert^2  \; \dx \geq 0.\] 
Now the identity $ \lambda \psi^{\prime \prime}(\lambda)= \psi^\prime(\lambda) $ for $\lambda\in \R$ implies 
\begin{align*}
    h_z^{\prime}(t)&\leq  \int_{\T^3}\psi^{ \prime}(\vr(t,x)-K) \lr{\vv-\mathbf{b}}\nabla(\vr(t,x)-K) \; \dx \\
    &=\int_{\T^3}\nabla(\psi(\vr(t,x)-K))  \lr{\vv-\mathbf{b}} \; \dx\\
    &=- \int_{\T^3} \psi(\vr(t,x)-K)  \Div \lr{\vv-\mathbf{b}} \; \dx \\
    &\leq \|\Div \lr{\vv-\mathbf{b}}\|_{L_\infty(\T^3)} h_z(t).
\end{align*}
Here, we apply Gr\"onwall's inequality along with $\inf_{x\in \T^3} \vr_0 =K$ to deduce
\begin{align}
    h_z(t)=0 \text{ for a.e. } t \in (0,T).
\end{align}
Therefore, we have $ \vr \geq K= \inf_{x\in \T^3} \vr_0 >0$ in $(0,T)\times \T^3 $.
\begin{flushright}
$\square$
\end{flushright}
Next, we state another lemma, related to the maximum and minimal principles of \eqref{CE:lin} with $b=0$. 
\begin{lem}\label{lem:mp-ld}
    Assume $\vr_0 \in H^2(\T^3)$,  $\vv \in L_\infty(0,T;H^2(\T^3))$ and $\va \in L_\infty(0,T;W^1_\infty) \cap L_2(0,T;H^2(\T^3))$ with $\va \geq c>0$. Then the unique solution of \eqref{CE:lin} with $\vb=0$ admits satisfying 
    \begin{align}\label{max-min}
			 \inf_{\T^3}\vr_0 \exp\left(-\int_0^t \Vert \Div \mathbf{v}(s) \Vert_{L_\infty(\T^3) } \; ds \right)
				& \leq  \vr(t,x) \leq  \sup_{\T^3}\vr_0  \exp\left(\int_0^t \Vert \Div \mathbf{v}(s) \Vert_{L_\infty(\T^3) } \; ds\right),
    \end{align}
    \text{ for } $0\leq t \leq T$. 
\end{lem}
{\emph{Sketch of Proof.}} The proof is  similar to the proof non-negativity property in Lemma \ref{lem:CE:exi} and extensively used in literatures like Novotn\'y and Straskraba \cite[Proposition 7.39]{NS}, Feireisl and Novotny \cite[Lemma 3.1]{FN}, where they have consider a special case with $\vc{a}$ is constant. For the sake of completeness we just highlight the key steps:  
			\begin{itemize}[leftmargin=*]
                \item Define : $R=(\sup_{\T^3}\vr_0)  \exp\left(\int_0^t \Vert \Div \mathbf{v}(s) \Vert_{L_\infty(\T^3) } \; ds\right).  $
				\item Then $R$ satisfies $R^{\prime}(t)- \Vert \Div \mathbf{v}(t) \Vert_{L_\infty(\T^3)} R(t) =0$ with $R(0)= \sup_{\T^3}\vr_0$ .
				\item Now consider $W(t,x)= \vr(t,x)- R(t)$ and it satisfies 
    \begin{align}\label{mp-1}
        \partial_t W + \Div(W\vv ) -\Div(\va \nabla \vr) \leq 0 \text{ a.e. in }(0,T) \times \R^d ,
    \end{align}
    with $W(0,x)= \vr_0- \sup_{\T^3}\vr_0\leq 0 $. 
    \item Now we test the equation \eqref{mp-1} with $\Psi^{\prime}(W)$ and integrating over space to obtain 
    \begin{align*}
        \frac{d}{dt} \int_{\T^3} \Psi(W) + 2 \int_{\T^3} \va  \Psi^{\prime \prime }(W) |\nabla W|^2 \leq \Vert \Div \mathbf{v}\Vert_{L_\infty(\T^3) } \Psi(W)
    \end{align*}
    where 
    \[\mathit{\Psi(\lambda)}= 
   \begin{cases}
   \frac{1}{2}\lambda^2, & \lambda\geq 0\\
   0, &\lambda< 0.\\
\end{cases}\]
\item Since, $\Psi(0)\leq 0$, from Gr\"onwall's inequality, we conclude 
\[  \vr(t,x) \leq  \sup_{\T^3}\vr_0  \exp\left(\int_0^t \Vert \Div \mathbf{v}(s) \Vert_{L_\infty(\T^3) } \; ds\right)\]  
\end{itemize}
For the other side of the inequality \eqref{max-min}, we need to consider 
\begin{align*}
    r(t)= \inf_{\T^3}\vr_0 \exp\left(-\int_0^t \Vert \Div \mathbf{v}(s) \Vert_{L_\infty(\T^3) } \; ds \right)
\end{align*}
and proceed analogously by considering 
\[ w(t,x)= \vr(t,x)-r(t) .\]

\begin{flushright}
$\square$
\end{flushright}
We have collected all the necessary tools to use in order to prove the convergence of iterative scheme \eqref{iter}. 


\section{Convergence of the iterative scheme}\label{conv}
\subsection{Boundedness of the sequence of approximations}
The estimates from the previous section allow
us to show 
\begin{lem}\label{lem:bd1}
Let $(\vw^n,\vr^n)$ be the sequence defined in \eqref{iter}
with $\vu^1(t,x) = \vu_0(x)$ for $(t,x) \in (0,T)\times \T^3$ and initial data satisfying \eqref{init1} or \eqref{init2}.  
There exists $M=M(\|\vr_0\|_{H^4(\T^3)},\|\vu_0\|_{H^3(\T^3)})$ in case of \eqref{init1} or $M=M(\|\vr_0\|_{H^3(\T^3)},\|\vu_0+\nabla p(\vr_0)\|_{H^3(\T^3)})$ in case of \eqref{init2} and $T=T(M)>0$ 
such that
\begin{equation} \label{est:seq}
\|\vw^n\|_{\CY_3(T)}+\|\vr^n\|_{\CV_3(T)}\leq M\quad \forall n \in \N,
\end{equation}
where the spaces $\CY_k(T)$ and $\CV_k(T)$ are defined in \eqref{def:CX}.
\end{lem}

\noindent{\em Proof.}  Recall that $\phi(\cdot)$ denotes an increasing, positive function, precise form of which may vary from line to line. For the purpose of this proof we also introduce more precise notation $\phi_1-\phi_3$ to denote given increasing, positive functions. The first equation of system \eqref{iter} is exactly \eqref{CE:lin} with $b=0$, $\vr=\vr^{n+1},\vv=\vw^{n+1}$ and $\va=\vr^n p'(\vr^n)$, while the subequation of the second line 
of \eqref{iter} corresponding to the $i$-th component of $\vw^{n+1}$ is nothing but \eqref{trans} with $g=0$, $\eta=\vw^{n+1}_i$ and $\vv=\vu^n$. 

By Lemma \ref{l:trans}, recalling that $\vu^n=\vw^n+\nabla p(\vr^n)$ we have
\begin{equation} \label{est:wn}
\|\vw^{n+1}\|_{\CY_3(T)}\leq \phi(\sqrt{T}\|\vu^n\|_{\CX_3(T)})=\phi_1\left(\sqrt{T}(\|\vw^n(t)\|_{\CY_3(T)}+\|\vr^n\|_{\CV_3(T)})\right)
\end{equation}
In order to use \eqref{est:CE} to estimate $\|\vr^{n+1}\|_{\CV_3(T)}$ we need to have a closer look at $\va(\vr)=\vr p'(\vr)$.  
\begin{equation*}
\begin{aligned}
&\nabla \va(\vr) = \left[p'(\vr)+\vr p^{(2)}(\vr)\right] \nabla \vr,\\ 
&\nabla^2 \va(\vr) = \left[p'(\vr)+\vr p^{(2)}(\vr)\right]\nabla^2\vr + \left[p^{(2)}(\vr)+\vr p^{(3)}(\vr)\right] |\nabla \vr|^2,\\
&\nabla^3 \va(\vr) = Q(\vr,p'(\vr),p^{(2)}(\vr),p^{(3)}(\vr),p^{(4)}(\vr)) [\nabla^3 \vr + \nabla^2 \vr \nabla \vr+|\nabla \vr|^3 ],  
\end{aligned}
\end{equation*}
where $Q$ is some polynomial, a precise form of which is not relevant. 
We have $\vr \in L_\infty( (0,T)\times \T^3)$. Therefore, as $p \in C^5$,  
\begin{equation*}
\begin{aligned}
&\|\nabla^2 \va(\vr) \|_{L_2(0,T;L_2(\T^3))}^2 \leq \phi(\|\vr\|_{\linfty})
\left[\|\nabla^2 \vr\|_{L_2(0,T;L_2(\T^3))}^2+\||\nabla \vr|^2\|_{L_2(0,T;L_2(\T^3))}^2\right] \\
&\leq 
T \phi(\|\vr\|_{\linfty}) \left(\|\nabla^2\vr\|_{L_\infty(0,T;H^1(\T^3))}^2 + C\|\nabla \vr\|_{\linfty}\right)
\end{aligned}
\end{equation*}
and 
\begin{equation*}
\begin{aligned}
&\|\nabla^3 \va(\vr)\|_{L_2(0,T;L_2(\T^3))}^2 \\
&\leq \phi(\|\vr\|_{\linfty})\left(
\|\nabla^3 \vr\|_{L_2(0,T;L_2(\T^3))}^2 + \|\nabla^2\vr \nabla\vr\|_{L_2(0,T;L_2(\T^3))}^2 + \||\nabla \vr|^3\|_{L_2(0,T;L_2(\T^3))}^2\right) \\
&\leq T \phi(\|\vr\|_{\linfty})\\ 
&\times \left(\|\nabla^3 \vr\|_{L_\infty(0,T;L_2(\T^3))}^2
+ \|\nabla\vr\|_{\linfty}^2\|\nabla^2 \vr\|_{L_\infty(0,T;L_2(\T^3))}^2
+ C \||\nabla\vr|^3\|_{\linfty} \right).
\end{aligned}
\end{equation*}
These bounds with an obvious estimate on $\|\nabla \va(\vr)\|_{L_2(0,T;L_2(\T^3))}$ imply 
\begin{equation} \label{a1}
\|\nabla \va(\vr)\|_{L_2(0,T;H^2(\T^3))} \leq T \phi(\|\vr\|_{\linfty}) \|\vr\|_{L_\infty(0,T;H^3(\T^3))}.
\end{equation}
Obviously we also have 
\begin{equation} \label{a2}
\|\nabla\va(\vr)\|_{\linfty}\leq \phi(\|\vr\|_{\linfty}) \|\vr\|_{L_\infty(0,T;H^3(\T^3))}.
\end{equation}
Applying \eqref{est:CE} to the first equation of \eqref{iter} we obtain  
\begin{equation*}
\begin{aligned}
&\|\vr^{n+1}\|_{\CV_3(T)} \leq \phi(\|\vr_0\|_{H^3(\T^3)}) \Big[C(T)+ \|\vr^{n+1}\|_{\CV_3(T)}\\ 
&\times[ T (\|\vw^{n+1}\|_{\CY_3(T)}
+ \|\nabla\va(\vr^n)\|_{\linfty}) +\|\nabla \va(\vr^n)\|_{L_2(0,T;H^2(\T^3))}]\Big],
\end{aligned}
\end{equation*} 
which together with \eqref{a1}, \eqref{a2} and \eqref{est:wn} gives
\begin{equation} \label{est:vrn}
\begin{aligned}
&\|\vr^{n+1}\|_{\CV_3(T)} \\ 
&\leq \phi_2(\|\vr_0\|_{H^3(\T^3)}) \left[C(T)+ 
T\|\vr^{n+1}\|_{\CV_3(T)}\left( \|\vw^{n+1}\|_{\CY_3(T)}+\phi_3(\|\vr^n\|_{\linfty})\|\vr^n\|_{\CV_3(T)} \right) \right]\\
&\leq \phi_2(\|\vr_0\|_{H^3(\T^3)}) \left[C(T)+ 
T\|\vr^{n+1}\|_{\CV_3(T)}\Big[ \phi_1\Big(\sqrt{T}(\|\vw^n(t)\|_{\CY_3(T)}+\|\vr^n\|_{\CV_3(T)})\right)\\
&\qquad \qquad \qquad \qquad+\phi_3(\|\vr^n\|_{\linfty})\|\vr^n\|_{\CV_3(T)} \Big]\Big]. 
\end{aligned}
\end{equation}
Let us take 
$$
M=2\max\{\sup_{s\in[0,1]}\phi_1(s),\phi_2(\|\vr_0\|_{H^3(\T^3)})C(T)\},
$$
where $C(T)$ is the constant from \eqref{est:vrn}.
Then, assuming that  
$$
\|\vw^n\|_{\CY_3(T)}+\|\vr^n\|_{\CV_3(T)}\leq M,
$$
for sufficiently small $T$ we can assure that 
\begin{align*}
&\phi_1(\sqrt{T}(\|\vw^n\|_{\CY_3(T)}+\|\vr^n\|_{\CV_3(T)})\leq\frac{M}{2}\\
& T\phi_2(\|\vr_0\|_{H^3})\left[ \|\vw^{n+1}\|_{\CY_3(T)}+\phi_3(\|\vr^n\|_{\linfty})\|\vr^n\|_{\CV_3(T)} \right] \leq \frac{1}{2},
\end{align*}
which together with \eqref{est:wn} and \eqref{est:vrn} implies 
$$
\|\vw^{n+1}\|_{\CY_3(T)}+\|\vr^{n+1}\|_{\CV_3(T)}\leq M,
$$
thus we have \eqref{est:seq}.
\begin{flushright}
$\square$
\end{flushright}

\subsection{Contraction argument. Proof of Theorem \ref{thm:main}. }\label{cont}
\begin{lem} \label{lem:contrac}
 Under the assumptions of Lemma \ref{lem:bd1} we have
\begin{equation} \label{contrac:final}
\|\delta \vw^{n+1}\|_{L_\infty(0,T;H^2(\T^3))}+\|\delta \vr^{n+1}\|_{\CV_2(T)}\leq E(T)\left( \|\delta \vw^n\|_{L_\infty(0,T;H^2(\T^3))}+\|\delta \vr^n\|_{\CV_2(T)} \right),
\end{equation}
where notation $E(t)$ is described in Section \ref{Notation}. 
\end{lem}
{\emph{Proof.}}
Let us denote 
\begin{equation*} 
\delta \vw^n = \vw^n-\vw^{n-1}, \quad \delta \vr^n = \vr^n-\vr^{n-1}, \quad 
\delta \vu^n = \vu^n-\vu^{n-1}.
\end{equation*}
Subtracting \eqref{iter} for $(\vw^{n+1},\vr^{n+1})$ and $(\vw^{n},\vr^{n})$ we obtain 
\begin{equation} \label{iter:dif}
\left\{ \begin{array}{lr}
\delta \vr^{n+1}_t + \div(\delta\vr^{n+1}\vw^{n+1})-\div(\vr^n p'( \vr^n)\nabla\delta\vr^{n+1}) = R_n ,\\
\delta\vw^{n+1}_t + \vu^n\cdot \nabla \delta\vw^{n+1} = -\delta \vu^n\cdot\nabla\vw^n,\\
(\delta\vr^{n+1},\delta\vw^{n+1})|_{t=0}= (0,\bf 0),
\end{array}\right.
\end{equation}
where
\begin{equation} \label{rn}
R_n = \div \left[ \big(p'(\vr^{n-1})\delta\vr^n + \vr^n(p'(\vr^n)-p'(\vr^{n-1}))\big)\nabla\vr^n -\vr^n\delta\vw^{n+1} \right].
\end{equation}
Each equation of the second line of \eqref{iter:dif} corresponds to \eqref{trans} with $g \sim \delta \vu^n \nabla \vw^n$. Therefore, taking into account \eqref{est:wn}, we can differentiate the right hand side in space only twice. For this purpose we show contraction in lower regularity then the estimate \eqref{est:wn}. This approach is well known in the regularity theory of the compressible and inhomogeneous Navier-Stokes systems to overcome the limitations coming from the presence of the gradient of the density in the continuity equation (see among others \cite{Hoff},\cite{MP},\cite{PP},\cite{DM}, \cite{DMP}). 
Combining \eqref{est:trans:01} and \eqref{est:trans:05} for $\eta_0=0$ we obtain
\begin{equation*}
\|\eta\|_{L_\infty(0,T;H^2(\T^3))} \leq E(T) \|g\|_{L_2(0,T;H^2(\T^3))},  
\end{equation*}
which applied to \eqref{iter:dif} implies 
\begin{equation} \label{est:dif1}
\|\delta \vw^{n+1}\|_{L_\infty(0,T;H^2(\T^3))} \leq E(T) \|\delta \vu^n\cdot\nabla\vw^n\|_{L_q(0,T;H^2(\T^3))} \quad \forall,\ 1<q\leq\infty.
\end{equation}

The first equation of \eqref{iter:dif} is \eqref{CE:lin} with $\vv=\vw^{n+1}$, $\va=\vr^n p'(\vr^n)$, $b=R_n$ and $\vr_0=0$.

Therefore \eqref{est:CE:0} implies 
\begin{equation} \label{est:dif2}
\begin{aligned}
&\|\delta\vr^{n+1}\|_{\CV_2(T)} \leq C \Big[\|R_n\|_{L_2(0,T;H^1(\T^3))}
+\|\delta\vr^{n+1}\|_{L_\infty(0,T;H^2(\T^3))}\\ 
&\times\Big( T \big(\|\vw^{n+1}\|_{L_\infty(0,T;H^2(\T^3))}
+ \|\nabla \big(\vr^n p'(\vr^n)\big) \|_{\linfty}\big) +\|\nabla^2 \big(\vr^n p'(\vr^n)\big)\|_{L_2(0,T;L_2(\T^3))}]\Big)\Big].
\end{aligned}
\end{equation}
As $p(\cdot)$ is sufficiently smooth, we have
$$
\nabla^2(\vr^np'(\vr^n)) \sim |\nabla \vr^n|^2+|\vr^n|\,|\nabla^2 \vr^n| 
$$
By \eqref{est:seq} we have
$$
\||\nabla\vr^n|^2\|_{L_2(0,T;L_2(\T^3))}^2=\int_0^T \|\nabla \vr^n\|_{L_4(\T^3)}^4 \leq T \|\nabla \vr^n\|_{L_\infty(0,T;L_4(\T^3))} 
\leq M\, E(T),
$$
and 
$$
\|\vr^n \nabla^2 \vr^n\|_{L_2(0,T;L_2(\T^3))}
\leq \|\vr^n\|_{\linfty}\sqrt{T}\|\nabla^2 \vr^n \|_{L_\infty(0,T;L_2(\T^3))}\leq M\, E(T). 
$$
Therefore choosing $T$ sufficiently small we can ensure that 
$$
\textrm{2nd line of \eqref{est:dif2}} \leq \frac{1}{2}. 
$$
Choosing such $T$ and adding \eqref{est:dif2} to \eqref{est:dif1} we obtain 
\begin{equation} \label{contrac:1}
\|\delta \vw^{n+1} \|_{L_\infty(0,T;W^2_2(\T^3))}+\|\delta \vr^{n+1}\|_{\CV_2(T)}\leq C \|R_n\|_{L_2(0,T;H^1(\T^3))}+E(T) \|\delta \vu^n\cdot\nabla\vw^n\|_{L_q(0,T;W^2_2(\T^3))} .
\end{equation}
Recalling \eqref{rn} we have 
\begin{align} \label{rn:2}
\|R_n\|_{L_2(0,T;H^1(\T^3))}\leq & \|p'(\vr^{n-1})\delta \vr^n \nabla \vr^n\|_{L_2(0,T;H^2(\T^3))}+\|\vr^n(p'(\vr^n)-p'(\vr^{n-1}))\nabla \vr^n\|_{L_2(0,T;H^2(\T^3))}\nonumber\\
&+\|\vr^n \delta\vw^{n+1}\|_{L_2(0,T;H^2(\T^3))}=:A_1+A_2+A_3.
\end{align}
The last term can be estimated directly: 
\begin{equation} \label{a3}
A_3 \leq \|\vr^n\|_{\linfty}\|\delta \vw^{n+1}\|_{L_2(0,T;H^2(\T^3))}\leq E(T)\|\delta \vw^{n+1}\|_{L_\infty(0,T;H^2(\T^3))},
\end{equation}
Let us proceed with $A_1$. We estimate the second order derivatives, which are the most restrictive. We have 
\begin{align*}
\nabla^2 \left( p'(\vr^n)\delta \vr^n \nabla \vr^n \right) \sim & \nabla^2 \vr^{n-1} \delta\vr^n \nabla\vr^n + \nabla \vr^{n-1} \delta \vr^n \nabla^2 \vr^n 
+ \nabla \vr^{n-1}\nabla \delta\vr^n \nabla \vr^n\\
&+\vr^{n-1}\nabla^2\delta \vr^n\nabla\vr^n 
+\vr^{n-1}\nabla\delta\vr^n\nabla^2 \vr^n+ \vr^{n-1}\delta \vr^n \nabla^3 \vr^n 
\end{align*}
The first two terms have the same structure and can be bounded as follows applying \eqref{est:seq}: 
\begin{align*}
\|\nabla^2 \vr^{n-1} \delta\vr^n \nabla\vr^n\|_{L_2(0,T;L_2(\T^3))} \leq & \|\nabla \vr^n\|_{\linfty} \|\nabla^2 \vr^{n-1}\|_{L_2(0,T;L_2(\T^3))}\|\delta \vr^n\|_{\linfty}\\
\leq & M^2E(T)\|\delta \vr^n\|_{\CV_2(T)}.
\end{align*}
For the third term we have  
\begin{align*}
\|\nabla \vr^{n-1}\nabla\delta\vr^n\nabla\vr^n\|_{L_2(0,T;L_2(\T^3))}\leq &
\|\nabla \vr^{n-1}\|_{\linfty}\|\nabla \vr^n\|_{\linfty}
\|\nabla \delta\vr^n\|_{L_2(0,T;L_2(\T^3))}\\\leq & M^2E(T)\|\delta \vr^n\|_{\CV_2(T)},
\end{align*}
and the same estimate holds for the fourth one. Next, 
\begin{align*}
\|\vr^{n-1}\nabla \delta\vr^n\nabla^2 \vr^n\|_{L_2(0,T;L_2(\T^3))} \leq & \sqrt{T} \|\vr^{n-1}\|_{\linfty}\|\nabla^2 \vr^n\|_{L_\infty(0,T;L_4(\T^3))}\|\nabla\delta \vr^n\|_{L_\infty(0,T;L_4(\T^3))}\\ 
\leq & M^2E(T)\|\delta \vr^n\|_{\CV_2(T)},
\end{align*}
and finally, for the last term, we have
\begin{align*}
\|\vr^{n-1}\delta \vr^n \nabla^3 \vr^n\|_{L_2(0,T;L_2(\T^3))}
& \leq C \|\vr^{n-1}\|_{\linfty} T \|\nabla^3 \vr^n\|_{L_\infty(0,T;L_2(\T^3))} \|\delta \vr^n\|_{\linfty}\\
& \leq M^2 E(T) \|\delta \vr^n\|_{\CV_2(T)}.
\end{align*}
Plugging the above estimates into \eqref{rn:2} we obtain and observing that $A_1$ and $A_2$ have the same structure due to assumed regularity of the pressure we obtain
\begin{equation*}
A_1+A_2 \leq M^2 E(T) \|\delta \vr^n\|_{\CV_2(T)},
\end{equation*}
which together with \eqref{a3} gives 
\begin{equation} \label{est:R1}
\|R_n\|_{L_2(0,T;H^1(\T^3))} \leq E(T) \left(\|\delta \vr^n\|_{\CV_2(T)} + \|\delta \vw^{n+1}\|_{L_\infty(0,T;H^2(\T^3))} \right).
\end{equation}
Now in order to close the contraction argument it is enough to estimate the second term on the RHS of \eqref{contrac:1}. Recalling \eqref{uw} we have 
\begin{equation*} 
\delta \vu^n = \delta \vw^n -  p'(\vr^n)\nabla\delta\vr^n - [p'(\vr^n)-p'(\vr^{n-1})]\nabla \vr^{n-1}. 
\end{equation*}
We have 
\begin{align}
\|\nabla^2 (\delta \vw^n \cdot \nabla \vw^n)\|_{L_2(0,T;L_2(\T^3))} & \leq 
\|\nabla^2 \delta\vw^n\cdot\nabla\vw^n\|_{L_2(0,T;L_2(\T^3))} + \|\nabla \delta\vw^n\cdot\nabla^2\vw^n\|_{L_2(0,T;L_2(\T^3))}\nonumber\\
& + \|\delta\vw^n\cdot\nabla^3\vw^n\|_{L_2(0,T;L_2(\T^3))} =: B_1+B_2+B_3  
\end{align}
We again apply \eqref{est:seq} to obtain 
\begin{align*}
&B_1 \leq \sqrt{T} \|\nabla \vw_n\|_{\linfty} \|\delta\vw^n\|_{L_\infty(0,T;H^2(\T^3))} \leq M E(T)\|\delta\vw^n\|_{L_\infty(0,T;H^2(\T^3))},\\
&B_2 \leq \sqrt{T} \|\vw^n\|_{L_\infty(0,T;H^2(\T^3))} \|\delta \vw^n\|_{L_\infty(0,T;H^2(\T^3))} \leq M E(T)\|\delta\vw^n\|_{L_\infty(0,T;H^2(\T^3))}, \\
&B_3 \leq \|\delta\vw^n\|_{\linfty} \|\nabla^3\vw^n\|_{L_2(0,T;L_2(\T^3))} \leq ME(T)\|\delta\vw^n\|_{L_\infty(0,T;H^2(\T^3))}, 
\end{align*}
which gives 
\begin{equation} \label{R2:1}
\|\delta \vw^n \cdot \nabla \vw^n\|_{L_2(0,T;H^2(\T^3))}
\leq E(T)\|\delta\vw^n\|_{L_\infty(0,T;H^2(\T^3))}.
\end{equation}
Next, 
\begin{align} \label{c2}
\nabla^2 \left( p'(\vr^n)\nabla \delta\vr^n \right)\sim
&\nabla^2 \vr^n \nabla \delta \vr^n \nabla \vw^n
+\nabla\vr^n\nabla^2\delta\vr^n+\nabla\vr^n\nabla\delta\vr^n\nabla^2\vw^n \nonumber\\
&+\vr^n\nabla^3\delta\vr^n\nabla\vw^n 
+\vr^n\nabla^2\delta\vr^n\nabla^2\vw^n
+\vr^n\nabla^3\delta\vr^n\nabla\vw^n.
\end{align}
Let us focus on the last term: 
\begin{align*}
    \|\vr^n \nabla^3 \delta \vr^n \cdot \nabla \vw^n\|_{L_2(0,T;L_2(\T^3))} & \leq C(T) \|\vr^n\|_{\linfty} \| \nabla \vw^n \|_{\linfty}  \| \nabla^3 \delta \vr^n \|_{L_2(0,T;L_2(\T^3))} \\
    &\leq CM^2 \| \delta \vr^n \|_{L_q(0,T;H^3(\T^3))}
\leq CM^2 \| \delta \vr^n \|_{\CV_2(T)} . 
\end{align*}

The lack of small constant is not a problem, since the term which we are estimating is already multiplied by a small constant in \eqref{contrac:1}. The other terms in 
\eqref{c2} can be estimated similarly, here we even get additional smallness in time. Summing up, we obtain 
\begin{equation} \label{R2:2}
\|p'(\vr^n)\nabla\delta\vr^n\cdot\nabla\vw^n\|_{L_2(0,T;H^2(\T^3))}\leq C(M,T)\|\delta\vr^n\|_{\CV_2(T)}
\end{equation}
Similarly we can show 
\begin{equation} \label{R2:3}
\big\|[p'(\vr^n)-p'(\vr^{n-1})]\nabla\vr^{n-1}\cdot\vw^n\big\|_{L_2(0,T;H^2(\T^3))}\leq E(T)\|\delta\vr^n\|_{\CV_2(T)}.
\end{equation}
Combining \eqref{R2:1},\eqref{R2:2} and \eqref{R2:3} we get 
\begin{equation} \label{est:R2}
\|\delta \vu^n\cdot\nabla\vw^n\|_{L_2(0,T;H^2(\T^3))}\leq E(T)\|\delta\vw^n\|_{L_\infty(0,T;H^2)}+C\|\delta\vr^n\|_{\CV_2(T)}.
\end{equation}
Plugging \eqref{est:R1} and \eqref{est:R2} into \eqref{contrac:1} we finally conclude \eqref{contrac:final}.
\begin{flushright}
$\square$
\end{flushright}

Now we complete the proof of Theorem \ref{thm:main} in a standard way. Inequality \eqref{contrac:final} implies  
$$
(\vw^n,\vr^n) \to (\vw,\vr) \;\; \textrm{strongly in} \; L_\infty(0,T;H^2(\T^3)) \times \CV_2(T).
$$
On the other hand, the estimate \eqref{est:seq} 
implies existence of a subsequence, which we can still denote $(\vw^n,\vr^n)$ 
$$
(\vw^n,\vr^n) \to (\vw,\vr) \;\; \textrm{weakly in} \; \CY_3(T) \times \CV_3(T).
$$
Setting $\vu=\vw-\nabla p(\vr)$ we easily verify that the limit satisfies \eqref{sys3}. 


\section{Existence theory for general velocity offsets} \label{general-offset}
The goal of this section is to extend the previous result to the case of different  velocity offset functions $p(\vr)$. 
\begin{itemize}[leftmargin=*]
    \item \textbf{Singular offset:} As mentioned in the introduction, Aceves et. al \cite{ABDM} derive a variant of dissipative Aw-Rascle system from a microscopic model of pedestrian dynamics. This corresponds to \eqref{sys3} with the offset function in the form \eqref{singoff}. Drawing motivation from this, we consider a more general offset function and the system \eqref{P:1}, namely 
\begin{align}\label{singoff2}
    p(\vr)= a \frac{\vr^\alpha}{(1-\vr)^\beta} \text{ with } a>0,\; \alpha>0 \text{ and }\; \beta >1. 
\end{align}
Similar form was considered in \cite{CNPZ} where the existence of regular solutions for certain approximation of this function was proven in one space dimension. The authors also performed the singular limit passage $a\to0$ obtaining in the limit the hard congestion system. Similar limit has been postulated in the multi-dimensional case, see \cite{ABDM}, but to our knowledge it has not yet been proved rigorously. Nevertheless, it is expected that in the hard congestion limit, in the saturated region where $\vr=1$, one cannot expect a regular solution due to appearance of extra forcing term. Therefore here we keep $a$ a positive constant, which results in restriction to unsaturated flow provided that the initial density is strictly below 1.
This case is discussed in Section \ref{Sec:4.1}.

\item \textbf{Non-local velocity offset:} It turns out that considering the  non-local offset function 
\begin{align}\label{nonlocoff}
    p(\vr)=K \ast \vr,
\end{align}
where $K$ is some non-local kernel, leads to reformulation of the system that generalizes the pressureless Euler-alignment model, see \cite{CPSZ2024} for further details.
Unfortunately, we are not able to apply our construction along with the linear theory (described in Sections 2 and 3) directly to the system \eqref{sys3} with \eqref{nonlocoff}. Instead, we consider a system where the offset function  $p(\vr)$ is a combination of local and a particular non-local component that corresponds to the Newtonian potential, i.e., 
\[ - \Delta K \ast \vr \approx \vr .\] 
More precisely, the closure relation is now of the form
\begin{align}
    \vw=\vu+  \Grad  p(\vr) + \Grad \Phi_\vr \text{ with } p \in C^5(\R_{+}) \text{ and }- \Delta \Phi_\vr = \vr -<\vr>,
\end{align}
and so \eqref{sys3} can be rewritten as follows
\begin{subnumcases}{\label{P-2}}
          \pt \vr+\Div(\vr \vw)-\Div(\vr \lr{ \Grad  p(\vr) + \Grad \Phi_\vr})=0,\label{P1-2}\\
	\pt(\vr \vw)+\Div(\vr \vw\otimes\vw)=\Div(\vr \vw\otimes \Grad p(\vr))+\Div(\vr \vw \otimes \nabla \Phi_\vr ),
	\label{P2-2}\\
	- \Delta \Phi_\vr = \vr -<\vr>.\label{P3-2}
\end{subnumcases}
This case is discussed in Section \ref{Sec:4.2}.
\end{itemize}

\subsection{Construction for the system with singular velocity offset}\label{Sec:4.1}
Here, we consider \eqref{P1:1}-\eqref{P2:1} with \eqref{singoff2}. Moreover, along with the hypothesis on initial data \eqref{init1} or \eqref{init2}, we need 
\begin{align}\label{ini_sing}
    0<\vr_0(x) < 1 \text{ for } x\in\T^3. 
\end{align}
  \begin{thm} \label{thm:2}
Assume the initial data satisfies \eqref{init1} or \eqref{init2} with \eqref{ini_sing}. Then there exists $T>0$ such that system \eqref{sys3} admits a unique solution $(\vr,\vw) \in \CV_3(T)\times \CY_3(T)$ with the estimate
$$
\|\vr\|_{\CV_3(T)}+\|\vw\|_{\CY_3(T)} \leq C(\|\vr_0\|_{H^4(\T^3)},\|\vu_0\|_{H^3(\T^3)}) \text{ and }  0<\vr < 1 \text{ in } [0,T)\times \T^3
$$
in case of \eqref{init1} or 
$$
\|\vr\|_{\CV_3(T)}+\|\vw\|_{\CY_3(T)} \leq C(\|\vr_0\|_{H^3(\T^3)},\|\vu_0+\nabla p(\vr_0)\|_{H^3(\T^3)}) \text{ and }  0<\vr < 1 \text{ in } [0,T)\times \T^3
$$
in case of \eqref{init2}.
\end{thm}
The proof of the Theorem \ref{thm:2} consists of two parts: the proof of boundedness of approximate solutions, and the compactness argument. The iteration scheme for this case is the same as in\eqref{iter}, but the proof of analogue of Lemma \eqref{lem:bd1} requires some alterations that we explain below. The contraction argument is similar to the one described in Section \ref{cont}, more precisely Lemma \eqref{lem:contrac}, and we skip the details here.
 

The main lemma for the uniform bounds reads:
 \begin{lem}\label{lem:bd2}
Let $(\vw^n,\vr^n)$ be the sequence defined in \eqref{iter}
with $\vu^1(t,x) = \vu_0(x)$ for $(t,x) \in (0,T)\times \T^3$.
There exists $M=M(\|\vr_0\|_{H^4(\T^3)},\|\vu_0\|_{H^3(\T^3)})$ in case of \eqref{init1} or $M=M(\|\vr_0\|_{H^3(\T^3)},\|\vu_0+\nabla p(\vr_0)\|_{H^3(\T^3)})$ in case of \eqref{init2} and $T(=T(M))>0$ 
such that
\begin{equation} \label{est:seq:2}
\|\vw^n\|_{\CY_3(T)}+\|\vr^n\|_{\CV_3(T)}\leq M \text{ and }  0<\vr_n(t,x) < 1 \; \; \forall n \in \N,\; \forall (t,x)\in (0,T)\times \T^3,
\end{equation}
where the spaces $\CY_k(T)$ and $\CV_k(T)$ are defined in \eqref{def:CX}. More precisely, it holds that there exists $0<\vt <1$ such that 
\begin{align}\label{rho-up-below}
     \frac{\vt}{2}<\vr_n(t,x) < 1-\frac{\vt}{2} \; \; \forall n \in \N,\; \forall (t,x)\in (0,T)\times \T^3. 
\end{align}
\end{lem}
{\emph{Proof.}} Similarly to  Lemma \eqref{lem:bd1}, we prove this lemma  with the help of an induction hypothesis:  we interpret the first equation in \eqref{iter} as equation \eqref{CE:lin}
with $\vb=0$, $\vr=\vr^{n+1},\vv=\vw^{n+1}$ and $\va=\vr^n p'(\vr^n)$; while the $i$-th row  of the second equation in \eqref{iter} as \eqref{trans} with $g=0$, $\eta=\vw^{n+1}_i$
and $\vv=\vu^n$. \par 
Clearly, there exists $0<\vt<1$ such that
\begin{align*}\label{ini-del}
    \vt \leq \inf_{\T^3} \vr_0 \quad \text{ and } \quad \sup_{\T^3} \vr_0 \leq 1-\vt.
\end{align*}
 Now, we assume that $(\vr_n, \vw_n)$ satisfy \eqref{est:seq:2} along with \eqref{rho-up-below}. This yields, 
\begin{align*}
    \va=\vr^n p'(\vr^n) >c(\vt) >0
\end{align*}
and $\va \in L_\infty(0,T;W^1_\infty(\T^3)) \cap L_2(0,T;H^2(\T^3))$. 
Here, we perform the similar estimates for $ \vw^{n+1}  $ and obtain 
\begin{equation} \label{est:wn:2}
\|\vw^{n+1}\|_{\CY_3(T)}\leq \phi(\sqrt{T}\|\vu^n\|_{\CX_3(T)})=\phi_1\left(\sqrt{T}(\|\vw^n(t)\|_{\CY_3(T)}+\|\vr^n\|_{\CV_3(T)})\right) .
\end{equation}
Similarly, for $ \vr^{n+1}  $ we have 
\begin{equation} \label{est:vrn:2}
\begin{aligned}
&\|\vr^{n+1}\|_{\CV_3(T)} \\ 
&\leq \phi_2(\|\vr_0\|_{H^3(\T^3)}) \left[C(T)+ 
T\|\vr^{n+1}\|_{\CV_3(T)}\Big[ \phi_1\Big(\sqrt{T}(\|\vw^n(t)\|_{\CY_3(T)}+\|\vr^n\|_{\CV_3(T)})\right)\\
&+\phi_3(\|\vr^n\|_{\linfty})\|\vr^n\|_{\CV_3(T)} \Big]\Big]. 
\end{aligned}
\end{equation}
Applying Lemma \ref{lem:mp-ld} we obtain
\begin{align*}
     \inf_{\T^3}\vr_0 \exp\left(-\int_0^t \Vert \Div \vw^{n+1}(s) \Vert_{L_\infty } \; ds \right)
				& \leq  \vr^{n+1}(t,x) \leq  \sup_{\T^3}\vr_0  \exp\left(\int_0^t \Vert \Div \vw^{n+1}(s) \Vert_{L_\infty } \; ds\right),
\end{align*}
for $t\in (0,T)$. 
Because
\begin{align*}
    \Vert \Div \vw^{n+1} \Vert_{L_\infty((0,T)\times \T^3) } \leq C \Vert  \vw^{n+1} \Vert_{ \CV_3(T)},
\end{align*}
with $C $ independent of time, we further deduce  
\begin{align}\label{bd2-1}
    \vt \exp\left(- t C \Vert  \vw^{n+1} \Vert_{ \CV_3(T)} \right)
				& \leq  \vr^{n+1}(t,x) \leq  (1-\vt) \exp\left(t C \Vert  \vw^{n+1} \Vert_{ \CV_3(T)}\right).
\end{align}
We set
$$
M=2\max\{\sup_{s\in[0,1]}\phi_1(s),\phi_2(\|\vr_0\|_{H^3(\T^3)})C(T)\},
$$
where $C(T)$ is the constant from \eqref{est:vrn:2}. Assuming that  
$$
\|\vw^n\|_{\CY_3(T)}+\|\vr^n\|_{\CV_3(T)}\leq M,
$$
for sufficiently small $T$, we can show that 
\begin{align*}
&\phi_1(\sqrt{T}(\|\vw^n\|_{\CY_3(T)}+\|\vr^n\|_{\CV_3(T)})\leq\frac{M}{2},\\
& T\phi_2(\|\vr_0\|_{H^3(\T^3)})\left[ \|\vw^{n+1}\|_{\CY_3(T)}+\phi_3(\|\vr^n\|_{\linfty})\|\vr^n\|_{\CV_3(T)} \right] \leq \frac{1}{2}. 
\end{align*}
This in turn, along with \eqref{est:wn:2} and \eqref{est:vrn:2}, implies 
$$
\|\vw^{n+1}\|_{\CY_3(T)}+\|\vr^{n+1}\|_{\CV_3(T)}\leq M,
$$
thus we have the first part of \eqref{est:seq:2}.
Finally, we choose a sufficiently small $T$ depending only on $M$, such that from \eqref{bd2-1} we obtain that 
\begin{align*}
    &\frac{\vt}{2}<\vr^{n+1}(t,x) < 1-\frac{\vt}{2}\quad  \text{ on } (0,T)\times \T^3. 
\end{align*}
This finishes the proof.
\begin{flushright}
$\square$
\end{flushright}

\subsection{Construction for the system with nonlocal velocity offset}\label{Sec:4.2}
In this part we prove the following existence result  for system \eqref{P-2}: 
 \begin{thm} \label{thm:3}
Assume the initial data satisfies \eqref{init1} or \eqref{init2}. Then there exists $T>0$ such that system \eqref{P-2} admits a unique solution $(\vr,\vw, \Phi_\vr) \in \CV_3(T)\times \CY_3(T) \times \CV_5(T)$ with the estimate
$$
\|\vr\|_{\CV_3(T)}+\|\vw\|_{\CY_3(T)} + \|\Phi_\vr\|_{\CV_5(T)} \leq C(\|\vr_0\|_{H^4(\T^3)},\|\vu_0\|_{H^3(\T^3)})
$$
in case of \eqref{init1} or 
$$
\|\vr\|_{\CV_3(T)}+\|\vw\|_{\CY_3(T)} + \|\Phi_\vr\|_{\CV_5(T)}\leq C(\|\vr_0\|_{H^3(\T^3)},\|\vu_0+\nabla p(\vr_0)\|_{H^3(\T^3)})
$$
in case of \eqref{init2}.
\end{thm} 

Again the proof of the above theorem is similar to the proof of Theorem \ref{thm:main}. 
    

In this case, however, we use  the following iteration scheme for  construction of the approximate solution: 
\begin{equation} \label{iter3}
\left\{ \begin{array}{lr}
\vr^{n+1}_t + \div(\vr^{n+1}\vw^{n+1})=\div(\vr^n p'( \vr^n)\nabla\vr^{n+1})+ \div(\vr^{n+1} \nabla\Phi_{n}),\\
\vw^{n+1}_t + \vu^n\cdot \nabla \vw^{n+1} = 0,\\
- \Delta \Phi_{n+1} = \vr^{n+1} -<\vr^{n+1}>,\\
(\vr^{n+1},\vw^{n+1})|_{t=0}= (\vr_0,\vu_0+\nabla p(\vr_0)+ \nabla \Phi_{\vr_0}),
\end{array}\right.
\end{equation}
where $\displaystyle - \Delta \Phi_{\vr_0} = \vr_0 -<\vr_0>$,  $(\vr^0,\vw^0,\Phi_0)=(\vr_0,\vu_0+\nabla p(\vr_0)+ \nabla \Phi_{\vr_0}, \Phi_{\vr_0})$ and 
\begin{align}\label{v-w-3}
    \vw^n=\vu^n+  \Grad  p(\vr^n) + \Grad \Phi_n.
\end{align}

The key observation ts that  the equation for $\Phi_{n+1}$ in this scheme is an elliptic equation. Therefore, the regularity class for $\{\Phi_n\}$ is higher.
Additionally, we have uniform bounds for the approximate sequence $ (\vr^{n}, \vw^n, \Phi_n)  $ in spaces $\CV_3(T)\times \CY_3(T) \times \CV_5(T)$, as stated in the following Lemma.
  
    
\begin{lem}\label{lem:bd3}
Let $(\vw^n,\vr^n,\Phi_n)$ be the sequence defined in \eqref{iter3}
with $\vu^1(t,x) = \vu_0(x)$ for $(t,x) \in (0,T)\times \T^3$.
There exists $M=M(\|\vr_0\|_{H^4(\T^3)},\|\vu_0\|_{H^3(\T^3)})$ in case of \eqref{init1} or $M=M(\|\vr_0\|_{H^3},\|\vu_0+\nabla p(\vr_0)\|_{H^3(\T^3)})$ in case of \eqref{init2} and $T=T(M)>0$
such that
\begin{equation*} \label{est:seq3}
\|\vw^n\|_{\CY_3(T)}+\|\vr^n\|_{\CV_3(T)}+\|\Phi_n \|_{\CV_5(T)} \leq M\quad \forall n \in \N,
\end{equation*}
where the spaces $\CY_k(T)$ and $\CV_k(T)$ are defined in \eqref{def:CX}.
\end{lem}
{\emph{Sketch of the proof.}} The  key steps of the proof are the following.
\begin{itemize}[leftmargin=*]
\item The bounds on $\|\vw^{n+1}\|_{\CY_3(T)}$ and $\|\vr^{n+1}\|_{\CV_3(T)}$ are obtained by calculations similar to proof of Lemma \ref{lem:bd1}, but keeping in mind that the relation between $\vw_n $ and $\vu_n$ has been modified \eqref{v-w-3}. Ultimately, we obtain
    \begin{equation} \label{est:wn:3}
\|\vw^{n+1}\|_{\CY_3(T)}\leq \phi(\sqrt{{T}}\|\vu^n\|_{\CX_3(T)})=\phi_1\left(\sqrt{{T}}(\|\vw^n(t)\|_{\CY_3(T)}+\|\vr^n\|_{\CV_3(T)} + \| \Phi_n \|_{\CV_5(T) } )\right) ,
\end{equation}
and 
\begin{equation} \label{est:vrn:3}
\begin{aligned}
&\|\vr^{n+1}\|_{\CV_3(T)} \\ 
&\leq \phi_2(\|\vr_0\|_{H^3(\T^3)}) \left[C(T)+ 
T\|\vr^{n+1}\|_{\CV_3(T)}\Big[ \phi_1\Big(\sqrt{{T}}(\|\vw^n(t)\|_{\CY_3(T)}+\|\vr^n\|_{\CV_3(T)} +\| \Phi_n \|_{\CV_5(T) } )\right)\\
&\hspace{6.5cm}+\phi_3(\|\vr^n\|_{\linfty})\lr{\|\vr^n\|_{\CV_3(T)}+\| \Phi_n \|_{\CV_5(T)}} \Big]\Big]. 
\end{aligned}
\end{equation}

    \item The estimate of $ \|\Phi_{n+1} \|_{\CV_5(T)} $ follows from the elliptic regularity estimates for the solutions of the problem
\[- \Delta \Phi_{n+1} = \vr^{n+1} -<\vr^{n+1}>.\]  
We obtain 
\begin{align}\label{est:phi:3}
    \| \Phi_{n+1} \|_{\CV_5(T) } \leq C_{\rm ell} \| \vr^{n+1} \|_{\CV_3(T) },
\end{align}
where $C_{\rm ell}$ is independent of $T$.

\item The strict positivity of $\vr^{n+1}$ is a direct consequence of Lemma \ref{lem:CE:exi}, where we use the particular form of $\vb = \Div(\vr^{n}  \nabla \Phi_{n} )$ in the first line of \eqref{iter3},
i.e.,
\begin{align*}
    \vr^{n+1}_t+\div(\vr^{n+1} \vw^{n+1} ) -\div(\vr^n p'( \vr^n)\nabla\vr^{n+1})=\vb= \div(\vr^{n+1} \nabla\Phi_{n} ).
\end{align*}
\item 
Now we combine estimates \eqref{est:wn:3}-\eqref{est:phi:3} and choose 
$$
M=2(1+C_{\rm ell}) \max\{\sup_{s\in[0,1]}\phi_1(s),\phi_2(\|\vr_0\|_{H^3(\T^3)})C(T)\},
$$
where $C(T)$ is the constant from \eqref{est:vrn:3} and $C_{\rm ell}$ is from \eqref{est:phi:3}.
Now, assuming that  
$$
\|\vw^n\|_{\CY_3(T)}+\|\vr^n\|_{\CV_3(T)}+ \|\Phi_{n} \|_{\CV_5(T)} \leq M,
$$
for sufficiently small $T$, first we can show that 
\begin{align*}
&\phi_1(\sqrt{T}(\|\vw^n\|_{\CY_3(T)}+\|\vr^n\|_{\CV_3(T)}+\|\Phi_{n+1} \|_{\CV_5(T)} )\leq\frac{M}{2}\\
& T\phi_2(\|\vr_0\|_{H^3(\T^3)})\left[ \|\vw^{n+1}\|_{\CY_3(T)}+\phi_3(\|\vr^n\|_{\linfty})\lr{\|\vr^n\|_{\CV_3(T)}+\|\Phi_{n} \|_{\CV_5(T)} } \right] \leq \frac{1}{2}. 
\end{align*}
This, along with \eqref{est:wn:2} and \eqref{est:vrn:2}, implies 
$$
\|\vw^{n+1}\|_{\CY_3(T)}+\|\vr^{n+1}\|_{\CV_3(T)}+\|\Phi_{n+1} \|_{\CV_5(T)}\leq M,
$$
and the proof is complete.
\begin{flushright}
$\square$
\end{flushright}
\end{itemize}

The second part of the proof of  Theorem \ref{thm:3}, is to provide an analogue of the contraction argument described in Section \ref{cont}. We denote 
\begin{equation*} \label{dw:3}
\delta \vw^n = \vw^n-\vw^{n-1}, \quad \delta \vr^n = \vr^n-\vr^{n-1}, \quad 
\delta \vu^n = \vu^n-\vu^{n-1}, \quad \delta \Phi_n= \Phi_n-\Phi_{n-1}.
\end{equation*}
The  system of difference equation corresponding to \eqref{iter3} is thus 
\begin{equation*} \label{iter:dif:3}
\left\{ \begin{array}{lr}
\delta \vr^{n+1}_t + \div(\delta\vr^{n+1} \vw^{n+1})-\div(\vr^n p'( \vr^n)\nabla\delta\vr^{n+1}) = \div(\delta\vr^{n+1} \Grad \Phi_{n})+ \widetilde{R}_n ,\\
\delta\vw^{n+1}_t + \vu^n\cdot \nabla \delta\vw^{n+1} = -\delta \vu^n\cdot\nabla\vw^n,\\
-\Delta \delta \Phi_{n+1}= \delta \vr^{n+1},\\
(\delta\vr^{n+1},\delta\vw^{n+1})|_{t=0}= (0,\bf 0),
\end{array}\right.
\end{equation*}
where
\begin{equation*}
\widetilde{R}_n = \div \left[ \big(p'(\vr^{n-1})\delta\vr^n + \vr^n(p'(\vr^n)-p'(\vr^{n-1}))\big)\nabla\vr^n -\vr^n\delta\vw^{n+1} \right] + \div(\vr^{n} \Grad \delta \Phi_{n}).
\end{equation*}
Performing the estimates similar to those from the proof of Lemma \ref{lem:contrac}, for the term  
\[\|\delta \vw^{n+1}\|_{L_\infty(0,T;H^2(\T^3))}+\|\delta \vr^{n+1}\|_{\CV_2(T)}+ \|\delta \Phi_{n+1} \|_{\CV_4(T)}.\]
we can show the following Lemma:
\begin{lem} \label{lem:contrac3}
 Under the assumptions of Lemma \ref{lem:bd3} we have
\begin{align*} 
\begin{split}
&\|\delta \vw^{n+1}\|_{L_\infty(0,T;H^2(\T^3))}+\|\delta \vr^{n+1}\|_{\CV_2(T)} + + \|\delta \Phi_{n+1} \|_{\CV_4(T)} \\
&\leq E(T)\left( \|\delta \vw^n\|_{L_\infty(0,T;H^2(\T^3))}+\|\delta \vr^n\|_{\CV_2(T)} ++ \|\delta \Phi_{n} \|_{\CV_4(T)}\right),
\end{split}
\end{align*}
where $E(t)$ is described in Section \ref{Notation}. 
\end{lem}
This concludes the contraction argument,  and also the proof of Theorem \eqref{thm:3}.
\begin{flushright}
$\square$
\end{flushright}

\medskip 
{\bf Acknowledgment:} The work of N.C. was funded by the ``Excellence Initiative Research University (IDUB)" program at the University of Warsaw. T.P.'s work was supported by the National Science Centre (NCN) project 2022/45/B/ST1/03432. The work of N.C. and E.Z. was supported by the EPSRC Early Career Fellowship no. EP/V000586/1.

\end{document}